\documentclass[pdflatex,sn-mathphys-num]{template}
\graphicspath{{figures/}}

\usepackage{graphicx}
\usepackage{multirow}
\usepackage{amsmath,amssymb,amsfonts}
\usepackage{amsthm}
\usepackage{mathrsfs}
\usepackage[title]{appendix}
\usepackage{xcolor}
\usepackage{textcomp}
\usepackage{manyfoot}
\usepackage{booktabs}
\usepackage{algorithm}
\usepackage{algorithmicx}
\usepackage{listings}
\usepackage{hyperref}
\usepackage{url}
\usepackage{booktabs}
\usepackage{amsfonts}
\usepackage{microtype}
\usepackage{color}
\usepackage{indentfirst}
\usepackage{epstopdf}
\usepackage{bm}
\usepackage{bbm}
\usepackage{mathtools}
\usepackage{latexsym,enumerate}
\usepackage{soul}
\usepackage{subfig}
\usepackage{adjustbox}
\usepackage{pifont}
\usepackage[english]{babel}
\usepackage{enumitem}
\usepackage{tabularx}
\usepackage{array}
\usepackage{makecell}
\usepackage{multirow,multicol}
\usepackage{caption}
\usepackage{algorithm}
\usepackage{algpseudocode}
\usepackage{fancybox}
\usepackage{cleveref}
\usepackage{subcaption}
\Crefname{theorem}{Theorem}{Theorems}
\Crefname{lemma}{Lemma}{Lemmas}
\Crefname{corollary}{Corollary}{Corollaries}
\Crefname{proposition}{Proposition}{Propositions}
\Crefname{assumption}{Assumption}{Assumptions}
\Crefname{definition}{Definition}{Definitions}
\Crefname{example}{Example}{Examples}
\Crefname{remark}{Remark}{Remarks}
\Crefname{figure}{Figure}{Figures}
\Crefname{table}{Table}{Tables}

\crefformat{equation}{(#2#1#3)}
\Crefformat{equation}{(#2#1#3)}
\crefrangeformat{equation}{(#3#1#4)--(#5#2#6)}
\Crefrangeformat{equation}{(#3#1#4)--(#5#2#6)}
\crefmultiformat{equation}{(#2#1#3)}{ and (#2#1#3)}{, (#2#1#3)}{, and (#2#1#3)}
\Crefmultiformat{equation}{(#2#1#3)}{ and (#2#1#3)}{, (#2#1#3)}{, and (#2#1#3)}

\theoremstyle{thmstyleone}
\newtheorem{theorem}{Theorem}[section]
\newtheorem{proposition}{Proposition}[section]
\newtheorem{lemma}{Lemma}[section]
\newtheorem{assumption}{Assumption}[section]
\theoremstyle{thmstyletwo}
\newtheorem{example}{Example}[section]
\newtheorem{remark}{Remark}[section]
\theoremstyle{thmstylethree}

\begin{document}

\title{Parameter-related strong convergence rates of Euler-type methods for time-changed stochastic differential equations}

\author*[1]{ \sur{Ruchun Zuo}}\email{ruchunzuo@outlook.com}

\affil*[1]{\orgdiv{Department of Mathematics}, \orgname{Shanghai Normal University}, \orgaddress{\city{Shanghai}, \postcode{200234}, \country{China}}}

\abstract{\unboldmath
	An Euler-type framework with equidistant step sizes is proposed for a class of time-changed stochastic differential equations.
	We establish the strong convergence rate of the standard Euler--Maruyama method under the global Lipschitz condition.
	The theoretical analysis is then extended to the truncated Euler--Maruyama method, proving its strong convergence under relaxed Khasminskii-type conditions.
	Under suitable conditions, the strong convergence orders of both numerical schemes are shown to be close to $\alpha/2$, where $\alpha \in (0,1)$ is the parameter of the time-change process.
	These results are significantly different from existing works using random step sizes, which typically preserve the classical convergence order of $1/2$.
	Numerical simulations are provided to demonstrate the theoretical findings.	}


\maketitle

\section{Introduction}\label{sec:intro}

Stochastic differential equations (SDEs) driven by time-changed Brownian motion have emerged as a powerful mathematical framework for modeling anomalous diffusion phenomena across diverse scientific disciplines. 
These models have found widespread applications in physics \cite{METZLER2001}, mathematical finance \cite{magdziarz2011option}, and biological systems \cite{saxton1997single}, where traditional diffusion models fail to capture the complex temporal patterns observed in experimental data. 
The key feature of time-changed SDEs is the presence of subordination processes, wherein a random time change governed by an inverse subordinator introduces memory effects and trapping events into the system dynamics \cite{kobayashi2011stochastic}. 
This non-Markovian structure allows for more realistic modeling of systems exhibiting subdiffusive behavior, characterized by slower-than-normal diffusion rates.
For the standard It\^o SDEs
\begin{equation*}
	\mathrm{d}X(t) = f(X(t)) \mathrm{d}t + g(X(t)) \mathrm{d}B(t),
\end{equation*}
the probability density function $p(x,t)$ of the solution $X(t)$ satisfies the Fokker--Planck equations (FPEs)
\begin{equation*}
	\frac{\partial p(x,t)}{\partial t} = \mathcal{L}_{\text{FP}} p(x,t),
\end{equation*}
where $\mathcal{L}_{\text{FP}} = -\frac{\partial}{\partial x} ( f(x) \cdot ) + \frac{1}{2} \frac{\partial^2}{\partial x^2} ( g(x)^2 \cdot )$. 
In contrast, when considering subdiffusion phenomena via time-changed SDEs, as rigorously established by \cite{gorenflo1997fractional, kobayashi2011stochastic}, the probability density function of a time-changed process $X(E(t))$ (to be formally defined subsequently) is governed by fractional Fokker--Planck equations (FFPEs)
\begin{equation*}
	D_t^\alpha p(x,t) = \mathcal{L}_{\text{FP}} p(x,t),
\end{equation*}
where $D_t^\alpha$ represents the Caputo fractional derivative of order $\alpha$ with respect to $t$. 
However, obtaining analytical solutions to the FFPEs presents significant mathematical challenges. 
Consequently, the solutions of the FFPEs can be characterized through the probability density function of solutions to the corresponding time-changed SDEs.

Given the analytical complexity of time-changed SDEs, numerical methods play a crucial role in their practical implementation. 
The first systematic investigation of numerical methods for time-changed SDEs was conducted by \cite{jum2014strong}, who considered the following time-changed SDEs
\begin{equation} \label{eq:dualitySDE1}
	\mathrm{d}Y(t) = f(E(t), Y(t)) \mathrm{d}E(t) + g(E(t), Y(t)) \mathrm{d}B(E(t)), \quad Y(0) = X_0,
\end{equation}
where $E(t)$ represents the inverse subordinator and $B(E(t))$ denotes the time-changed Brownian motion. By establishing a duality principle that connects the time-changed SDEs to the classical It\^o SDEs
\begin{equation} \label{eq:dualitySDE2}
	\mathrm{d}X(t) = f(t, X(t)) \mathrm{d}t + g(t, X(t)) \mathrm{d}B(t), \quad X(0) = X_0,
\end{equation}
they showed that, under global Lipschitz conditions, if $X(t)$ solves \Cref{eq:dualitySDE2}, then $Y(t) = X(E(t))$ solves \Cref{eq:dualitySDE1}, and conversely, if $Y(t)$ solves \Cref{eq:dualitySDE1}, then $X(t) = Y(D(t))$ solves \Cref{eq:dualitySDE2}. 
Through this duality relationship, they developed an Euler--Maruyama (EM) method achieving a strong convergence rate of $1/2$.

Following this pioneering work, \cite{jin2019strong} investigated a class of time-changed SDEs where the inverse subordinator $E(t)$ in the drift and diffusion coefficients of \Cref{eq:dualitySDE1} is replaced with ordinary time $t$. 
In this modified framework, the duality principle no longer applies, necessitating the development of alternative analytical techniques. 
Using a Gronwall-like inequality specific to time-changed integrals, they demonstrated that the EM method maintains a strong convergence rate of $1/2$ in this setting as well.

Subsequently, numerical methods for time-changed SDEs have been extensively studied under different discretization strategies. 
These include backward EM methods \cite{deng2020semi}, truncated EM methods \cite{liu2020truncated, li2023convergence, li2025truncated}, and higher-order methods such as the Milstein scheme \cite{jin2021strong, liu2023milstein}. 
A notable commonality among these approaches is their reliance on random discretization of the inverse subordinator $E(t)$, which effectively transforms the stochastic increments of the time-change process into deterministic increments within the numerical scheme. 
As a consequence, these methods generally yield convergence rates similar to those established for classical SDEs (typically of rate $1/2$ or $1$), despite the additional complexity introduced by the time-change process.

A natural question arises: can we develop numerical schemes for time-changed SDEs whose convergence properties explicitly reflect the unique characteristics of the time-change process? 
	This question motivates the central innovation of our present work: developing equidistant-step Euler-type methods for time-changed SDEs. 
	Unlike previous approaches, our methods preserve the intrinsic stochasticity of the subordinator process increments within the numerical scheme, leading to convergence behavior that directly depends on the stability index $\alpha$ of the subordinator.

It is worth noting that an equidistant-step EM method has been previously utilized by \cite{carnaffan2017solving} in their investigation of anomalous diffusion processes, which provide stochastic representations for multidimensional FFPEs.
Their approach employed techniques from Malliavin calculus to derive unbiased density formulas for the probability density functions of these anomalous diffusion processes.
This methodology yielded remarkably accurate numerical approximations that successfully captured essential characteristics of the underlying processes, including the distinctive cusp formations observed in density profiles.
Despite these promising applications, a rigorous mathematical analysis concerning the convergence properties of this numerical scheme in the context of time-changed SDEs has not been thoroughly established in the literature.
The theoretical results established in our present work thus provide a rigorous mathematical foundation for such numerical applications.

In this paper, we consider the following class of time-changed SDEs
\begin{equation}\label{eq:SDEs}
	\mathrm{d}X(t) = f(X(t)) \mathrm{d}E(t) + g(X(t)) \mathrm{d}B(E(t)),
\end{equation}
where $f: \mathbb{R}^d \rightarrow \mathbb{R}^d$, $g: \mathbb{R}^d \rightarrow \mathbb{R}^{d \times m}$, $E(t)$ denotes the inverse subordinator, and $B(E(t))$ represents the $m$-dimensional time-changed Brownian motion.

The main contribution of this paper is the development and rigorous analysis of equidistant-step Euler-type methods for the time-changed SDEs \Cref{eq:SDEs}.
We first propose the standard equidistant-step EM method and establish its strong convergence under the global Lipschitz condition. 
However, such an assumption is often too restrictive for practical applications. 
In particular, when the drift and diffusion coefficients exhibit super-linear growth, the standard EM method is known to diverge in general \cite{hutzenthaler2011divergence}. 
To overcome this fundamental limitation, we relax the theoretical assumptions to Khasminskii-type conditions and further formulate an equidistant-step truncated EM scheme.
Under suitable conditions, the convergence orders of both numerical schemes are arbitrarily close to $\alpha/2$.
These explicit error estimates characterize how the stability index $\alpha$ of the subordinator process determines the convergence behavior, thereby revealing a fundamental relationship between the statistical properties of the time-change process and the accuracy of the numerical methods.

The remainder of this paper is organized as follows: 
\Cref{sec:pre} introduces the mathematical preliminaries, including key properties of subordinators and their inverses, and establishes important lemmas used in our subsequent analysis. 
\Cref{sec:EM} develops the equidistant-step EM scheme and proves its strong convergence rate. 
\Cref{sec:TEM} introduces the equidistant-step truncated EM scheme and establishes its strong convergence rate.
\Cref{sec:num} presents numerical experiments that validate our theoretical findings and illustrate the effectiveness of the proposed methods.

\section{Preliminaries}\label{sec:pre}

Throughout this paper, we work within a complete probability space $(\Omega, \mathcal{F}, \mathbb{P})$. 
We use $|\cdot|$ to denote the Euclidean norm of a vector and the Frobenius norm of a matrix.
For $a,b\in\mathbb{R}$, we write $a\wedge b=\min\{a,b\}$ and $a\vee b=\max\{a,b\}$.
We next introduce the fundamental stochastic processes that underlie our analysis of time-changed SDEs.

Let $D = (D(t))_{t \geq 0}$ denote a non-decreasing L\'evy process with $D(0) = 0$, characterized by c\`adl\`ag (right-continuous with left limits) sample paths. 
Such a process $D$ is called a subordinator and is uniquely determined by its Laplace exponent $\psi: [0,\infty) \to [0,\infty)$ through the relation
\begin{equation*}
	\mathbb{E}[e^{-sD(t)}] = e^{-t\psi(s)}, \quad s \geq 0, \quad t \geq 0.
\end{equation*}
The Laplace exponent $\psi$, with zero killing rate, is represented via the L\'evy--Khintchine formula as
\begin{equation*}
	\psi(s) = \int_{0}^{\infty} (1 - e^{-sy})   \nu(\mathrm{d}y), \quad s > 0,
\end{equation*}
where $\nu$ is the L\'evy measure satisfying the integrability condition $\int_{0}^{\infty} (y \wedge 1)   \nu(\mathrm{d}y) < \infty$. 
Throughout this paper, we assume that $\nu(0, \infty) = \infty$, which ensures that $D$ has strictly increasing paths with infinitely many jumps.
The inverse subordinator $E = (E(t))_{t \geq 0}$ associated with $D$ is defined by
\begin{equation*}
	E(t) := \inf \{ u > 0 : D(u) > t \}, \quad t \geq 0.
\end{equation*}
This process $E$ is continuous, non-decreasing, and initialized at $E(0) = 0$. 
The fundamental relationship between $D$ and $E$ is characterized by the equivalence
\begin{equation*}
	\{E(t) > x\} = \{D(x) < t\}, \quad \forall t, x \geq 0.
\end{equation*}
In this work, we focus specifically on the case where $D$ is an $\alpha$-stable subordinator with stability index $\alpha \in (0,1)$, for which the Laplace exponent takes the simplified form $\psi(s) = s^{\alpha}$. 
Consequently, $E$ is referred to as the inverse $\alpha$-stable subordinator. 
The assumption $\nu(0, \infty) = \infty$ ensures that $E$ has continuous sample paths. 
A critical feature of this process is that jumps in $D$ correspond to intervals of constancy in $E$, during which any time-changed process $X \circ E = (X(E(t)))_{t \geq 0}$ remains constant.

\begin{figure}[htbp]
	\centering
	\includegraphics[width=0.48\linewidth]{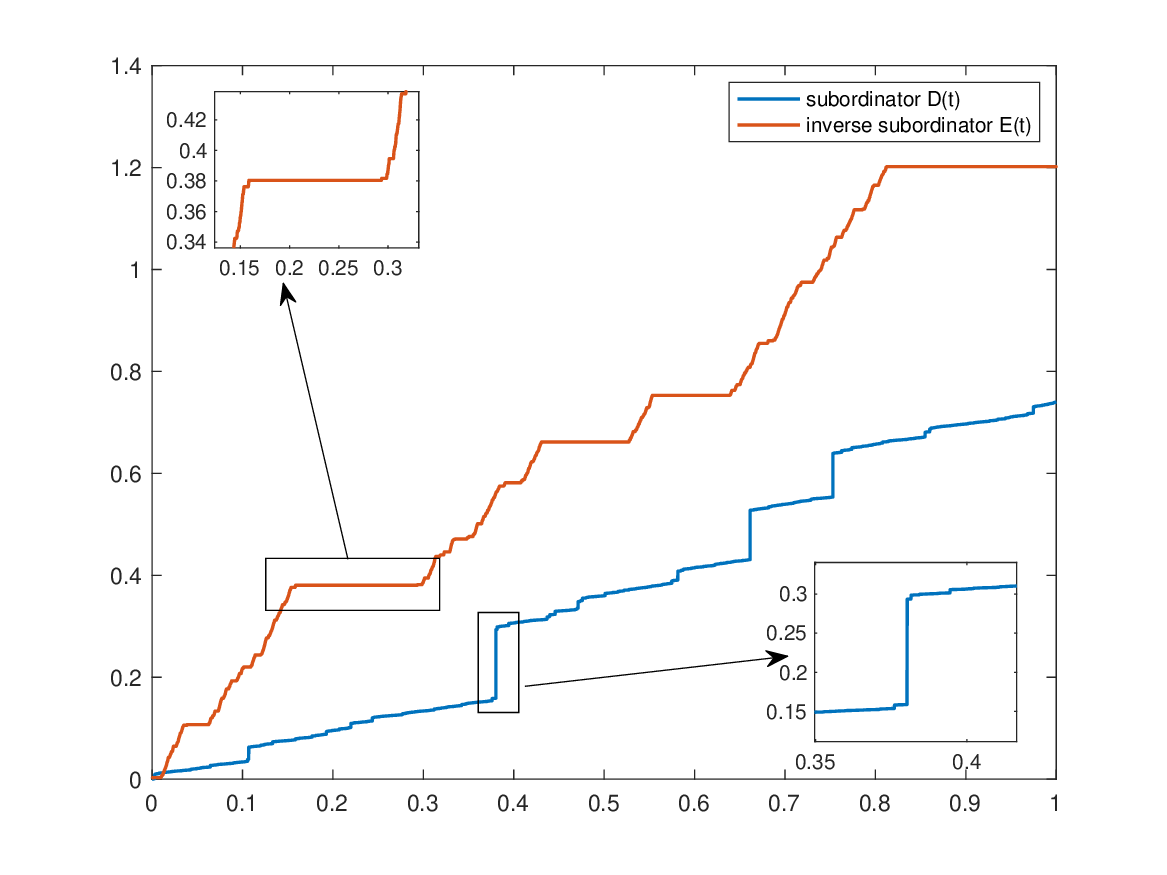}
	\hfill
	\includegraphics[width=0.48\linewidth]{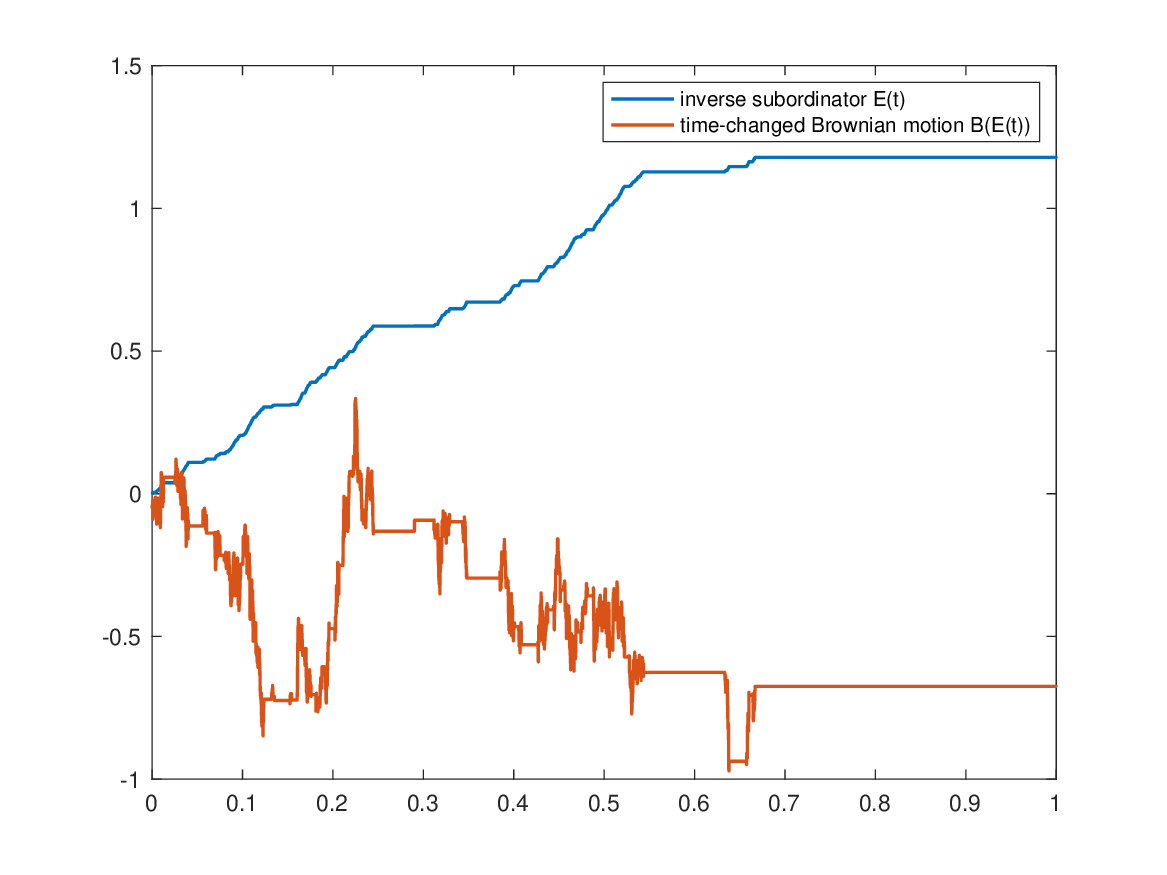}
	\caption{The left figure is sample paths of an 0.8-stable subordinator $D$ (blue) and its inverse subordinator $E$ (red). The right figure is sample paths of an inverse 0.8-stable subordinator $E$ (blue) and the corresponding time-changed Brownian motion $B \circ E$ (red).}
	\label{fig:D_E_BE}
\end{figure}

Let $B = (B(t))_{t \geq 0}$ be a standard $m$-dimensional Brownian motion defined on $(\Omega, \mathcal{F}, \mathbb{P})$ and independent of the subordinator $D$. 
We denote by $\mathbb{P}_B$ and $\mathbb{P}_D$ the probability measures associated with $B$ and $D$, respectively, and by $\mathbb{E}_B$ and $\mathbb{E}_D$ the corresponding expectation operators. 
The total expectation $\mathbb{E}$ under the probability measure $\mathbb{P}$ can then be expressed as
\begin{equation*}
	\mathbb{E}[\cdot] = \mathbb{E}_D[\mathbb{E}_B[\cdot]] = \mathbb{E}_B[\mathbb{E}_D[\cdot]].
\end{equation*}
The time-changed process $B \circ E = (B(E(t)))_{t \geq 0}$ represents a Brownian motion evaluated along the random time scale defined by $E$. 
This construction effectively models particles that remain stationary during the constant intervals of $E$, capturing the characteristic trapping events observed in subdiffusive systems. 
It is important to note that while $B \circ D$ retains the Brownian motion property from $B$, the time-changed process $B \circ E$ is non-Markovian, reflecting the memory effects introduced by the inverse subordinator.
The quadratic variations of these processes satisfy the following relationships
\begin{equation*}
	[B \circ E, B \circ E]_t = E(t), \quad [B \circ E, E]_t = [E, B \circ E]_t = 0, \quad t \geq 0.
\end{equation*}
For a comprehensive treatment of stochastic calculus for time-changed semimartingales, we refer the reader to \cite{kobayashi2011stochastic}. 
\Cref{fig:D_E_BE} provides a visual representation of the processes $D$, $E$, and $B \circ E$.

To numerically simulate the inverse subordinator $E$, we employ the method developed by \cite{magdziarz2009langevin, magdziarz2009stochastic}. 
Given a time horizon $T > 0$ and an equidistant step size $\delta > 0$, we first generate a sample path of $D$ over $[0, T]$ by setting $D(0) = 0$ and recursively computing
\begin{equation*}
	D(i\delta) = D((i-1)\delta) + Z_i, \quad i = 1, 2, \ldots,
\end{equation*}
where $(Z_i)_{i \in \mathbb{N}}$ constitutes an independent and identically distributed sequence with $Z_i \sim \Delta D_{\delta}$. 
These increments are generated using the formula
\begin{equation*}
	Z_i = \frac{\sin(\alpha(V + \pi/2))}{(\cos(V))^{1/\alpha}} \Bigl(\frac{\cos(V - \alpha(V + \pi/2))}{W}\Bigr)^{(1-\alpha)/\alpha},
\end{equation*}
where $V \sim \text{Uniform}(-\pi/2, \pi/2)$ and $W \sim \text{Exponential}(1)$ are independent random variables.
The approximation $\tilde{E} = (\tilde{E}(t))_{t \geq 0}$ of the inverse subordinator is then defined by
\begin{equation*}
	\tilde{E}(t) := ( \min \{ n \in \mathbb{N} : D(n\delta) > t \} - 1 ) \delta, \quad t \geq 0.
\end{equation*}
This approximation $\tilde{E}$ is a non-decreasing step function with jump size $\delta$. 
As demonstrated in \cite{magdziarz2009langevin, magdziarz2009stochastic}, it provides a close approximation to $E$ with the error bounds
\begin{equation*}
	E(t) - \delta \leq \tilde{E}(t) \leq E(t), \quad \text{for all } t \in [0, T],
\end{equation*}
holding almost surely.

We now present two fundamental lemmas from \cite{jum2014strong,jin2019strong,jin2021strong} regarding the moments and exponential moments of the inverse subordinator, which play a key role in our convergence analysis.

\begin{lemma}\label{lem:E_moment}
	Let $E$ denote the inverse of the $\alpha$-stable subordinator $D$.
	Then, for each $t \ge 0$ and $n \in \mathbb{N}$,
	\begin{equation*}
		\mathbb{E}[E^n(t)]= \frac{n!}{\Gamma(n\alpha+1)} t^{n\alpha},
	\end{equation*}
	where $\Gamma(\cdot)$ is the gamma function.
\end{lemma}

\begin{lemma}\label{lem:ExpRV}
	Let $E$ denote the inverse $\alpha$-stable subordinator with index $\alpha \in (0, 1)$.
	For any $t > 0$, $C > 0$, and $p > 0$, if $p < (1 - \alpha)^{-1}$, then it holds that
	\begin{equation*}
		\mathbb{E}[e^{C E(t)^p}] < \infty.
	\end{equation*}
	\end{lemma}

A critical challenge in the error analysis of equidistant-step numerical methods for time-changed SDEs arises from the fact that the inverse subordinator $E$ lacks independent and stationary increments \cite{meerschaert2004limit}.
Consequently, standard pointwise moment bounds for $E(t+h) - E(t)$ established in \cite{gajda2014large} are insufficient for estimating integrals driven by $E$, as the increments cannot be decoupled from the integrator.
To overcome this essential difficulty, we construct a uniform bound on the maximum increment of $E$ over the entire time horizon.
The following lemma establishes an upper bound on the moments of this supremum increment.
This uniform estimate constitutes a cornerstone of our strong convergence proof, allowing us to extract the maximal increment outside the integrals and effectively relate the approximation error to the time step size.

\begin{lemma}\label{lem:sup_delta_E}
	Let $E$ be the inverse $\alpha$-stable subordinator with index $\alpha\in(0,1)$. 
	For any $T>0$ and $h>0$, define the maximum increment over intervals of length $h$ by
	\begin{equation*}
		\Xi_{h}:=\sup_{0\le u\le T- h}(E(u+h)-E(u)).
	\end{equation*}
	Then, for any $n\ge 1$ and $\varepsilon\in(0,\alpha)$, there exists a constant $C>0$, independent of $h$, such that for all $h\in(0, 1/2]$,
	\begin{equation*}
		\mathbb{E}[\Xi_{h}^n]\le C(h)^{n(\alpha-\varepsilon)}.
	\end{equation*}
\end{lemma}
	
	\begin{proof}
		We first recall the uniform modulus of continuity for the inverse stable subordinator. 
		Consider the random variable defined by
		\begin{equation*}
			M := \sup_{0\le v<u\le 1/2}\frac{E(u)-E(v)}{(u-v)^\alpha|\log(u-v)|^{1-\alpha}}.
		\end{equation*}
		According to \cite{iksanov2017fractionally}, $M$ is almost surely finite and possesses finite exponential moments; in particular, $\mathbb{E}[M^n] < \infty$ for all $n \ge 1$.
		
		Restricting the supremum to intervals of length $h \in (0, 1/2]$, we have almost surely,
		\begin{equation*}
			\sup_{0\le t\le 1/2-h} \big(E(t+h)-E(t)\big) \le M h^\alpha |\log h|^{1-\alpha}.
		\end{equation*}
		Taking the $n$-th moment yields
		\begin{equation}\label{eq:half_interval_bound}
			\mathbb{E}\Bigl[\sup_{0\le t\le 1/2-h}(E(t+h)-E(t))^n\Bigr]\le \mathbb{E}[M^n]h^{n\alpha}|\log h|^{n(1-\alpha)}.
		\end{equation}
		
		To extend this estimate to an arbitrary interval $[0, T]$, we invoke the self-similarity property of the inverse stable subordinator: $\{E(ct)\}_{t\ge 0} \stackrel{d}{=} \{c^\alpha E(t)\}_{t\ge 0}$ for any $c>0$, where $\stackrel{d}{=}$ denotes equality in distribution. 
		By choosing $c=2T$ and setting $\eta = h/(2T)$, we obtain the distributional equality
		\begin{equation*}
			\Xi_{h}	\stackrel{d}{=} (2T)^\alpha \sup_{0\le s\le 1/2-\eta}\big(E(s+\eta)-E(s)\big).
		\end{equation*}
		Note that for $h \in (0, T]$, we have $\eta \in (0, 1/2]$. 
		Combining this scaling relation with \Cref{eq:half_interval_bound}, we deduce
		\begin{align*}
			\mathbb{E}\bigl[\Xi_{h}^n\bigr]
			&= (2T)^{n\alpha}\mathbb{E}\Bigl[\sup_{0\le s\le 1/2-\eta}(E(s+\eta)-E(s))^n\Bigr] \\
			&\le C\eta^{n\alpha}|\log\eta|^{n(1-\alpha)} \\
			&= Ch^{n\alpha}\Bigl|\log\frac{h}{2T}\Bigr|^{n(1-\alpha)}
		\end{align*}
		For any fixed $\varepsilon \in (0, \alpha)$ and sufficiently small $h$, the logarithmic term is dominated by the polynomial decay, i.e., $|\log(h/2T)|^{n(1-\alpha)} \le C h^{-n\varepsilon}$. 
		Consequently,
		\begin{equation*}
			\mathbb{E}\bigl[\Xi_{h}^n\bigr]
			\le C   h^{n(\alpha-\varepsilon)},
		\end{equation*}
		where $C$ is a positive constant independent of $h$. 
		This concludes the proof.
\end{proof}

\section{Strong convergence of EM method}\label{sec:EM}

In this section, we formulate the equidistant-step EM scheme for the time-changed SDEs \Cref{eq:SDEs} and rigorously establish its strong convergence rate. 
To proceed with our theoretical analysis, we first impose the following standard assumptions on the drift and diffusion coefficients.
\begin{assumption}\label{assum:lipschitz}
	There exists a constant $K>0$ such that  for all $x,y \in \mathbb{R}^d$,
	\begin{equation*}
		|f(x) - f(y)|  + |g(x) - g(y)| \le K|x - y|.
	\end{equation*}
\end{assumption}

\begin{assumption}\label{assum:linear_growth}
	There exists a constant $K>0$ such that for all $x \in \mathbb{R}^d$,
	\begin{equation*}
		|f(x)| + |g(x)| \le K(1 + |x|).
	\end{equation*}
\end{assumption}

Under these assumptions, the time-changed SDEs \cref{eq:SDEs} admits a unique strong solution \cite{kobayashi2011stochastic}.
We now formulate the equidistant-step EM method. 
Given a uniform partition $\{t_n\}_{n=0}^N$ of $[0, T]$ with step size $\Delta t = t_{n+1} - t_n$, the EM scheme is defined recursively by
\begin{equation}\label{eq:EM_discrete}
	X_{n+1} = X_n + f(X_n) \Delta E_n + g(X_n) \Delta B_n, \quad n = 0, 1, \ldots, N-1,
\end{equation}
where $\Delta E_n = E(t_{n+1}) - E(t_n)$ and $\Delta B_n = B(E(t_{n+1})) - B(E(t_n))$ represent the increments of the inverse subordinator and time-changed Brownian motion, respectively.
To facilitate theoretical analysis, we introduce a piecewise constant interpolation of the numerical solution, defined by
\begin{equation*}
	\tilde{X}_t = X_n \quad \text{for } t \in [t_n, t_{n+1}), \quad n = 0, 1, \ldots, N-1.
\end{equation*}
For any $s \in [0, T]$, let $n_s$ denote the index such that $s \in [t_{n_s}, t_{n_s+1})$. 
The discrete scheme \Cref{eq:EM_discrete} can then be expressed in the equivalent integral form
\begin{equation}\label{eq:EM_continuous}
	\tilde{X}_{t_{n_s}} = X_0 + \int_0^{t_{n_s}} f(\tilde{X}_r)   \mathrm{d}E(r) + \int_0^{t_{n_s}} g(\tilde{X}_r)   \mathrm{d}B(E(r)).
\end{equation}
Furthermore, we introduce a continuous-time interpolation of the numerical solution, defined by
	\begin{equation}\label{eq:inter_X}
		\bar{X}_t = X_0 + \int_0^{t} f(\tilde{X}_r)   \mathrm{d}E(r) + \int_0^{t} g(\tilde{X}_r)   \mathrm{d}B(E(r)).
\end{equation}

The following lemma, adapted from \cite{jin2019strong}, establishes the moment boundedness of the true solution to the time-changed SDEs \Cref{eq:SDEs}.
\begin{proposition}\label{pro:true_solution_p_bound}
	Let $ X $ be the solution of the time-changed SDEs \Cref{eq:SDEs}, where $ f $ and $g$ satisfy \Cref{assum:lipschitz} and \Cref{assum:linear_growth}. 
	Then, for any $ p \ge 1 $, we have $ \mathbb{E}[R_p(T)] < \infty $, where $ R_p(t) = 1 + \sup_{0 \le s \le t} |X(s)|^p $.
\end{proposition}

Next, we establish uniform moment bounds for the continuous-time interpolation \Cref{eq:inter_X}.
Before proceeding, we specify that throughout the following estimates, $C$ represents a positive constant independent of $\Delta t$ that may change from line to line.
\begin{lemma}\label{lem:int_X_monment}
	Suppose \Cref{assum:linear_growth} holds. For any $T>0$ and $p\ge 2$, there exists a constant $C>0$, independent of $\Delta t$, such that for all $\Delta t \in (0, 1]$,
	\begin{equation*}
		\mathbb{E}\Bigl[\sup_{0\le t\le T}|\bar{X}_t|^p\Bigr] \le C.
	\end{equation*}
\end{lemma}
	
\begin{proof}
	Fix $t \in [0, T]$. 
	Applying It\^o's formula {\color{red}\cite{umarov2018beyond}} to the function $|x|^p$, we obtain
	\begin{equation}\label{eq:Xbar_ito}
		|\bar{X}_t|^p = |X_0|^p + \int_0^t \mathcal{L}\bar{X}_s   \mathrm{d}E(s) + p \int_0^t |\bar{X}_s|^{p-2}\bar{X}_s^\top g(\tilde{X}_s)   \mathrm{d}B(E(s)),
	\end{equation}
	where $\mathcal{L}\bar{X}_s := p|\bar{X}_s|^{p-2}\bar{X}_s^\top f(\tilde{X}_s) + \frac{1}{2}p(p-1)|\bar{X}_s|^{p-2}|g(\tilde{X}_s)|^2$.
	
	Taking the supremum over $[0, T]$ and applying $\mathbb{E}_B$, we estimate the drift term using \Cref{assum:linear_growth} and Young's inequality,
	\begin{equation}\label{eq:Xbar_drift}
		\mathbb{E}_B\Bigl[\sup_{0\le t\le T}\int_0^t \mathcal{L}\bar{X}_s   \mathrm{d}E(s)\Bigr] \le C \mathbb{E}_B\Bigl[\int_0^T \bigl(1 + |\bar{X}_s|^p + |\tilde{X}_s|^p\bigr) \mathrm{d}E(s)\Bigr].
	\end{equation}
	For the stochastic integral term, the Burkholder--Davis--Gundy (BDG) inequality {\color{red}\cite{umarov2018beyond}} implies
	\begin{align*}
		I_{mart} &:= \mathbb{E}_B\Bigl[\sup_{0\le t\le T}\Bigl|\int_0^t p|\bar{X}_s|^{p-2}\bar{X}_s^\top g(\tilde{X}_s)\mathrm{d}B(E(s))\Bigr|\Bigr] \\
		&\le C   \mathbb{E}_B\Bigl[\Bigl(\int_0^T |\bar{X}_s|^{2p-2}|g(\tilde{X}_s)|^2 \mathrm{d}E(s)\Bigr)^{1/2}\Bigr].
	\end{align*}
	Using \Cref{assum:linear_growth} and Young's inequality, we have
	\begin{align}\label{eq:Xbar_diffusion}
		I_{mart} &\le C \mathbb{E}_B\Bigl[\bigl(\sup_{0\le u\le T}|\bar{X}_u|^p\bigr)^{\frac12}\Bigl(\int_0^T   |\bar{X}_s|^{p-2}|g(\tilde{X}_s)|^2 \mathrm{d}E(s)\Bigr)^{1/2} \Bigr] \nonumber\\
		&\le \frac{1}{2} \mathbb{E}_B\Bigl[\sup_{0\le u\le T}|\bar{X}_u|^p\Bigr] + C \mathbb{E}_B\Bigl[\int_0^T \bigl(1 + |\bar{X}_s|^p + |\tilde{X}_s|^p\bigr) \mathrm{d}E(s)\Bigr].
	\end{align}
	Substituting \Cref{eq:Xbar_drift,eq:Xbar_diffusion} into \Cref{eq:Xbar_ito},
	\begin{equation*}
		\mathbb{E}_B\Bigl[\sup_{0\le u\le T}|\bar{X}_u|^p\Bigr]\le C\bigl(1+|X_0|^p\bigr)+C\int_0^T \bigl(1+\mathbb{E}_B\bigl[\sup_{0\le u\le s}|\bar{X}_u|^p\bigr]\bigr)\mathrm{d}E(s),
	\end{equation*}
	where we used the fact that $|\tilde{X}_s|^p \le \sup_{0\le u\le s}|\bar{X}_u|^p$.
	
	By applying Gronwall's inequality, we obtain
	\begin{equation*}
		\mathbb{E}_B\Bigl[\sup_{0\le u\le T}|\bar{X}_u|^p\Bigr]	\le C \bigl(1 + |X_0|^p + E(T)\bigr) e^{C E(T)}.
	\end{equation*}
	Finally, taking $\mathbb{E}_D$ on both sides, and using \Cref{lem:E_moment,lem:ExpRV}, we conclude that
	\begin{equation*}
		\mathbb{E}\Bigl[\sup_{0\le t\le T}|\bar{X}_t|^p\Bigr] \le C \mathbb{E}\bigl[(1+E(T))e^{C E(T)}\bigr] < \infty.
	\end{equation*}
	This completes the proof.
\end{proof}

By the definition $\tilde{X}_t=\bar{X}_{t_{n_t}}$ for $t\in[0,T]$, it immediately follows that
\begin{equation*}
	\sup_{0<\Delta t\le1}\mathbb{E}\Bigl[\sup_{0\le t\le T}|\tilde{X}_t|^p\Bigr]	\le	\sup_{0<\Delta t\le1}\mathbb{E}\Bigl[\sup_{0\le t\le T}|\bar{X}_t|^p\Bigr]	\le C.
\end{equation*}

We now turn our attention to the strong convergence rate of \Cref{eq:EM_continuous}.
Before proceeding, we recall the notion of strong convergence for numerical schemes. 
A numerical approximation $X_{\Delta t}$ of a stochastic process $X$ is said to converge strongly with order $\gamma > 0$ if there exists a constant $C$, independent of the step size $\Delta t$, such that
\begin{equation*}
	\mathbb{E}\Bigl[\sup_{t \in [0,T]} |X(t) - X_{\Delta t}(t)| \Bigr] \leq C(\Delta t)^{\gamma},
\end{equation*}
where $T$ is the terminal time.

\begin{theorem}\label{th:discrete_EM_convergence}
	Suppose \Cref{assum:lipschitz,assum:linear_growth} hold.
	Assume further that $\alpha \in (1/2,1)$ and $2 \le p < (1-\alpha)^{-1}$.
	Then, for any fixed $\varepsilon \in (0, \alpha)$,  there exists a constant $C > 0$ independent of $\Delta t$ such that
	\begin{equation*}
		\mathbb{E}\Bigl[ \sup\limits_{t \in [0,T]}|X({t}) - \tilde{X}_{t}|^p\Bigr] \leq C\Delta t^{\frac{p(\alpha-\varepsilon)}{2}}.
	\end{equation*}
\end{theorem}

\begin{proof}
	To measure the deviation between the true solution of the SDEs and the numerical solution, we introduce the error term
	\begin{equation*}
		Z(t) := \sup_{s \in [0, t]} |X(s) - \tilde{X}_{t_{n_s}}|,
	\end{equation*}
	where $X(s)$ represents the true solution at time $s$ and $\tilde{X}_{t_{n_s}}$ is the numerical solution at time $t_{n_s}$.
	
	By comparing \Cref{eq:SDEs} and \Cref{eq:EM_continuous}, we obtain
	\begin{align*}
		X(s) - \tilde{X}_{t_{n_s}} &= \int_{0}^{t_{n_s}}  (f(X(r)) - f(\tilde{X}_r))  \mathrm{d}E(r) + \int_{0}^{t_{n_s}}  (g(X(r)) - g(\tilde{X}_r))    \mathrm{d}B(E(r)) \\
		&\quad + \int_{t_{n_s}}^{s} f(X(r))   \mathrm{d}E(r) + \int_{t_{n_s}}^{s} g(X(r))   \mathrm{d}B(E(r)),
	\end{align*}
	where $ t_{n_s} $ denotes the discrete time point such that $ s \in [t_{n_s}, t_{n_s+1}) $.
	From this decomposition, we derive the following bound for the error
	\begin{equation*}
		Z(t) \leq I_1 + I_2 + I_3 + I_4, \quad Z^p(t) \leq 4^{p-1}(I^p_1 + I^p_2 + I^p_3 + I^p_4),
	\end{equation*}
	where the terms are defined as
	\begin{align*}
		I_1 &= \sup_{s \in [0, t]} \Bigl|\int_{0}^{t_{n_s}} (f(X(r)) - f(\tilde{X}_r))    \mathrm{d}E(r) \Bigr|, \\
		I_2 &= \sup_{s \in [0, t]} \Bigl|\int_{0}^{t_{n_s}} (g(X(r)) - g(\tilde{X}_r))    \mathrm{d}B(E(r)) \Bigr|, \\
		I_3 &= \sup_{s \in [0, t]} \Bigl|\int_{t_{n_s}}^{s}  f(X(r))    \mathrm{d}E(r) \Bigr|, \\
		I_4 &= \sup_{s \in [0, t]} \Bigl|\int_{t_{n_s}}^{s}  g(X(r))    \mathrm{d}B(E(r)) \Bigr|.
	\end{align*}
	
	For $I_1$, by \Cref{assum:lipschitz}, we have 
	\begin{align*}
		I_1 &\leq \sup_{s \in [0, t]} \int_{0}^{t_{n_s}} \big|f(X(r)) - f(\tilde{X}_r) \big|   \mathrm{d}E(r) \\
		&\leq \sup_{s \in [0, t]} \int_{0}^{t_{n_s}} K|X(r) - \tilde{X}_r|   \mathrm{d}E(r)\\
		& \leq K\int_{0}^{t} Z(r)   \mathrm{d}E(r).
	\end{align*}
	Since $t_{n_t} \leq t \leq T$, using H\"older's inequality and taking the expectation $\mathbb{E}_B$ on both sides, we have
	\begin{align}\label{eq:I1}
		\mathbb{E}_B[I_1^{p}] &\leq K^{p}E^{p-1}(T) \int_{0}^{t} \mathbb{E}_B[Z^{p}(r)]   \mathrm{d}E(r).		
	\end{align}
	
	For $I_2$, by \Cref{assum:lipschitz}, the BDG inequality and H\"older's inequality, we have
	\begin{align}\label{eq:I2}
		\mathbb{E}_B[I_2^{p}] 
		&\leq \mathbb{E}_B\Bigl[\Bigl(\int_{0}^{t_{n_t}} \big|g(X(r)) - g(\tilde{X}_r) \big|^2   \mathrm{d}E(r)\Bigr)^{\frac{p}{2}}\Bigr] \nonumber\\
		&\leq \mathbb{E}_B\Bigl[\Bigl(\int_{0}^{t_{n_t}} K^2 |X(r) - \tilde{X}_r|^2   \mathrm{d}E(r)\Bigr)^{\frac{p}{2}}\Bigr] \nonumber\\
		&\leq K^{p}E^{\frac{p}{2}}(T) \int_{0}^{t} \mathbb{E}_B[Z^{p}(r)]   \mathrm{d}E(r).
	\end{align}
	
	For $I_3$, recalling the definition of $\Xi_{\Delta t}$ from \Cref{lem:sup_delta_E} and $R_p(T)$ from \Cref{pro:true_solution_p_bound}, and applying \Cref{assum:linear_growth}, we have
	\begin{align*}
		I_3 &\leq \sup_{s \in [0, t]} \int_{t_{n_s}}^{s} |f(X(r))|   \mathrm{d}E(r) \\
		&\le K R_1(T) \sup_{s\in[0,t]}\big(E(s)-E(t_{n_s})\big) \\
		&\le K R_1(T)   \Xi_{\Delta t}.
	\end{align*}
	Taking the $p$-th moment leads to
	\begin{equation}\label{eq:I3}
		\mathbb{E}_B[I_3^p] \le K^p   \mathbb{E}_B[R_p(T)]   \Xi_{\Delta t}^p.
		\end{equation}
	
	For $I_4$, we use the second change-of-variable formula for stochastic integrals \cite{umarov2018beyond}.
	For $s \in [0,t]$, the stochastic integral can be rewritten as
	\begin{equation*}
		\int_{t_{n_s}}^{s}g(X(r))\mathrm{d}B(E(r))
		= \int_{E(t_{n_s})}^{E(s)}g\bigl(X(D(v-))\bigr)\mathrm{d}B(v),
	\end{equation*}
	where $D$ is the subordinator associated with $E$. 
	Define the continuous martingale
	\begin{equation*}
		M(u) := \int_{0}^{u}g\bigl(X(D(v-))\bigr)\mathrm{d}B(v), \qquad 0\le u\le E(T).
	\end{equation*}
	Then $I_4(t)$ can be bounded by the modulus of continuity of $M$
	\begin{equation*}
		I_4 = \sup_{0\le s\le t}\bigl|M(E(s))-M(E(t_{n_s}))\bigr|
		\le \sup_{\substack{0\le r<u\le E(T) \\ u-r\le \Xi_{\Delta t}}}|M(u)-M(r)|.
	\end{equation*}
	Invoking the modulus of continuity estimate for continuous martingales (\cite[Theorem 7]{jin2019strong}), we have
	\begin{align}\label{eq:I4}
		\mathbb{E}_B[I_4^p]
		&\le C_p   \mathbb{E}_B\bigl[R_p(T)\bigr]  
		\Xi_{\Delta t}^{\frac p2} \Bigl(\log \frac{4E(T)}{\Xi_{\Delta t}}\Bigr)^{\frac p2}.
	\end{align}
	
	Combining \Crefrange{eq:I1}{eq:I4}, we obtain
	\begin{align*}
			\mathbb{E}_B[Z^p(t)]
			&\le C \Bigl(E^{p-1}(T) + E^{\frac p2-1}(T)\Bigr) \int_0^t \mathbb{E}_B[Z^p(r)]\mathrm{d}E(r) \\
			&\quad + C  \mathbb{E}_B[R_p(T)] \Bigl(\Xi_{\Delta t}^{p} + \Xi_{\Delta t}^{\frac p2}\Bigl(\log \frac{4E(T)}{\Xi_{\Delta t}}\Bigr)^{\frac p2}\Bigr).
	\end{align*}
	Applying Gronwall's inequality, we get
	\begin{equation*}
		\mathbb{E}_B[Z^p(T)] \le C   \mathbb{E}_B[R_p(T)] \Bigl(\Xi_{\Delta t}^{p} + \Xi_{\Delta t}^{\frac p2}\Bigl(\log \frac{4E(T)}{\Xi_{\Delta t}}\Bigr)^{\frac p2}\Bigr) e^{C(E^{p}(T)+E^{\frac p2}(T))}.
	\end{equation*}
	To simplify the notation, let
	\begin{equation*}
		A := \mathbb{E}_B[R_p(T)]   e^{C(E^p(T)+E^{\frac p2}(T))}, \qquad
		L_{\Delta t} := \log \frac{4E(T)}{\Xi_{\Delta t}}.
	\end{equation*}
	Taking $\mathbb{E}_D$ on both sides, we have
	\begin{equation*}
		\mathbb{E}[Z^p(T)] \le C   \mathbb{E}\bigl[A   \Xi_{\Delta t}^{p}\bigr] + C   \mathbb{E}\bigl[A   \Xi_{\Delta t}^{\frac p2} L_{\Delta t}^{\frac p2}\bigr].
	\end{equation*}
	Applying the Cauchy--Schwarz inequality yields
	\begin{align*}
		\mathbb{E}\bigl[Z^p(T)\bigr]
		&\le C\bigl(\mathbb{E}[A^2]\bigr)^{\frac12}\bigl(\mathbb{E}\bigl[\Xi_{\Delta t}^{2p}\bigr]\bigr)^{\frac12}
		+ C\bigl(\mathbb{E}[A^2]\bigr)^{\frac12}\bigl(\mathbb{E}\bigl[\Xi_{\Delta t}^{p}L_{\Delta t}^{p}\bigr]\bigr)^{\frac12}.
	\end{align*}
	Using \Cref{lem:ExpRV} and \Cref{pro:true_solution_p_bound}, combined with the assumption $p < (1-\alpha)^{-1}$, we have $\mathbb{E}[A^2] \le C$.
	Consequently,
	\begin{equation*}
		\mathbb{E}\bigl[Z^p(T)\bigr] \le C\biggl[\bigl(\mathbb{E}\bigl[\Xi_{\Delta t}^{2p}\bigr]\bigr)^{\frac{1}{2}}+\Bigl(\mathbb{E}\Bigl[\Xi_{\Delta t}^{p}\bigl(\log\frac{4E(T)}{\Xi_{\Delta t}}\bigr)^{p}\Bigr]\Bigr)^{\frac{1}{2}}\biggr].
	\end{equation*}
	By \Cref{lem:sup_delta_E}, for any $\varepsilon \in (0,\alpha)$, we know that 
	\begin{equation*}
		(\mathbb{E}[\Xi_{\Delta t}^{2p}])^{1/2} \le C \Delta t^{p(\alpha-\varepsilon)}.
	\end{equation*}
	For the second term, we use the property that for any $\eta \in (0,1)$, $\log x \le C_\eta x^\eta$ for $x \ge 1$. This implies
	\begin{equation*}
		\Xi_{\Delta t}^{p}\Bigl(\log\frac{4E(T)}{\Xi_{\Delta t}}\Bigr)^{p}\le CE^{\eta p}(T)\Xi_{\Delta t}^{(1-\eta)p}.
	\end{equation*}
	Applying the Cauchy--Schwarz inequality, \Cref{lem:E_moment,lem:sup_delta_E} (with $\varepsilon/2$ in place of $\varepsilon$), we have
	\begin{align*}
		\Bigl(\mathbb{E}\Bigl[\Xi_{\Delta t}^{p}\bigl(\log\frac{4E(T)}{\Xi_{\Delta t}}\bigr)^{p}\Bigr]\Bigr)^{\frac{1}{2}}
		&\le C\bigl(\mathbb{E}[E^{2\eta p}(T)]\bigr)^{\frac{1}{4}}\bigl(\mathbb{E}\bigl[\Xi_{\Delta t}^{2(1-\eta)p}\bigr]\bigr)^{\frac{1}{4}} \\
		&\le C\Delta t^{\frac{p(1-\eta)(\alpha-\varepsilon/2)}{2}}.
	\end{align*}
	By choosing $\eta := \varepsilon/(2\alpha)$, the exponent becomes $\frac{p}{2}(\alpha-\varepsilon+\frac{\varepsilon^2}{4\alpha})$. 
	Since $\frac{\varepsilon^2}{4\alpha} > 0$ and $\Delta t \le 1$, it follows that$$C\Delta t^{\frac{p(1-\eta)(\alpha-\varepsilon/2)}{2}} \le C\Delta t^{\frac{p(\alpha-\varepsilon)}{2}}.$$
	Thus, we conclude
	\begin{equation*}
		\mathbb{E}\bigl[Z^p(T)\bigr]
		\le C\Bigl(\Delta t^{p(\alpha-\varepsilon)}+\Delta t^{\frac{p(\alpha-\varepsilon)}{2}}\Bigr)
		\le C\Delta t^{\frac{p(\alpha-\varepsilon)}{2}}.
	\end{equation*}
	This completes the proof.
\end{proof}

Before proving the strong convergence of the continuous-time interpolation \Cref{eq:inter_X}, we first establish a uniform bound on the difference between the piecewise constant approximation \Cref{eq:EM_continuous} and its continuous counterpart \Cref{eq:inter_X}.

\begin{lemma}\label{lem:discrete_minus_continuous}
	Suppose \Cref{assum:lipschitz,assum:linear_growth} hold. 
	For any $p \ge 2$ and any fixed $\varepsilon \in (0, \alpha)$, there exists a constant $C > 0$, independent of $\Delta t$, such that for all sufficiently small $\Delta t > 0$,
	\begin{equation*}
		\mathbb{E}\Bigl[\sup_{0\le t\le T}|\bar{X}_t-\tilde{X}_t|^p\Bigr]
		\le C\Delta t^{\frac{p(\alpha-\varepsilon)}{2}}.
	\end{equation*}
\end{lemma}

\begin{proof}
	Recall the definition of $\Xi_{\Delta t}$ from \Cref{lem:sup_delta_E}. 
	For any $t \in [0, T]$, let $t_n$ be the grid point such that $t \in [t_n, t_{n+1})$. 
	Since $\tilde{X}_t = \bar{X}_{t_n}$ on this interval, the difference process satisfies
	\begin{equation*}
		\bar{X}_t - \tilde{X}_t = f(\bar{X}_{t_n})\big(E(t)-E(t_n)\big) + g(\bar{X}_{t_n})\big(B(E(t))-B(E(t_n))\big).
	\end{equation*}
	By \Cref{assum:linear_growth}, we have
	\begin{equation*}
		|f(\bar{X}_{t_n})| \vee |g(\bar{X}_{t_n})|\le C\Bigl(1+\sup_{0\le s\le T}|\bar{X}_s|\Bigr).
	\end{equation*}
	Consequently,
	\begin{equation*}
		\sup_{0\le t\le T}|\bar{X}_t-\tilde{X}_t|\le C\Bigl(1+\sup_{0\le s\le T}|\bar{X}_s|\Bigr)(J_1+J_2),
	\end{equation*}
	where the increment terms are defined by
	\begin{equation*}
		J_1 := \sup_{0\le t\le T} \big(E(t) - E(t_{n_t})\big) \le \Xi_{\Delta t},
	\end{equation*}
	and
	\begin{equation*}
		J_2 := \sup_{0\le t\le T} \big|B(E(t)) - B(E(t_{n_t}))\big|.
	\end{equation*}
	Taking the $p$-th moment and applying the Cauchy--Schwarz inequality, we obtain
	\begin{align*}
		\mathbb{E}\Bigl[\sup_{0\le t\le T}|\bar{X}_t-\tilde{X}_t|^p\Bigr]
		&\le C\mathbb{E}\Bigl[\Bigl(1+\sup_{0\le s\le T}|\bar{X}_s|^p\Bigr)(J_1^p+J_2^p)\Bigr] \\
		&\le C\Bigl(\mathbb{E}\Bigl[1+\sup_{0\le s\le T}|\bar{X}_s|^{2p}\Bigr]\Bigr)^{1/2}\bigl(\mathbb{E}\bigl[J_1^{2p}+J_2^{2p}\bigr]\bigr)^{1/2}.
	\end{align*}
	By \Cref{lem:int_X_monment}, the moment of the numerical solution is bounded uniformly in $\Delta t$. 
	Thus, it suffices to estimate the moments of the increments.
	
	For $J_1$, we directly use \Cref{lem:sup_delta_E} to get
	\begin{equation*}
		\mathbb{E}[J_1^{2p}] \le \mathbb{E}[\Xi_{\Delta t}^{2p}] \le C \Delta t^{p(\alpha-\varepsilon)}.
	\end{equation*}
	For $J_2$, conditioned on the subordinator $E$, the term represents the modulus of continuity of a Brownian motion over intervals of length at most $\Xi_{\Delta t}$. 
	Following the same reasoning as in the estimation of $I_4$ in \Cref{th:discrete_EM_convergence}, we have
	\begin{equation*}
		\mathbb{E}\bigl[J_2^{2p}\bigr] \le C\mathbb{E}\Bigl[\Xi_{\Delta t}^{p}\Bigl|\log\frac{4E(T)}{\Xi_{\Delta t}}\Bigr|^{p}\Bigr].
	\end{equation*}
	As established in the proof of \Cref{th:discrete_EM_convergence}, the logarithmic term is absorbed by the polynomial decay for any small $\varepsilon > 0$, yielding
	\begin{equation*}
		\bigl(\mathbb{E}[J_2^{2p}]\bigr)^{1/2} \le C\Delta t^{\frac{p(\alpha-\varepsilon)}{2}}.
	\end{equation*}
	Combining these estimates completes the proof.
\end{proof}

By combining the error bound for the discrete approximation in \Cref{th:discrete_EM_convergence} and the continuous interpolation estimate in \Cref{lem:discrete_minus_continuous}, we now proceed to establish our main strong convergence theorem.
\begin{theorem}\label{th:cotinuous_EM_convergence}
	Suppose \Cref{assum:lipschitz,assum:linear_growth} hold. 
	Let $\alpha \in (1/2, 1)$ and assume $2 \le p < (1-\alpha)^{-1}$. 
	Then, for any fixed $\varepsilon \in (0, \alpha)$, there exists a constant $C > 0$, independent of $\Delta t$, such that for all sufficiently small $\Delta t > 0$,
	\begin{equation*}
		\mathbb{E}\Bigl[\sup_{0\le t\le T}|X(t)-\bar{X}_t|^p\Bigr]
		\le C\Delta t^{\frac{p(\alpha-\varepsilon)}{2}}.
	\end{equation*}
\end{theorem}

\begin{remark}
	The restriction $\alpha>1/2$ ensures that the interval $[2,(1-\alpha)^{-1})$ is nonempty.
	The upper bound $p<(1-\alpha)^{-1}$ is used in the proof of \Cref{th:discrete_EM_convergence} when applying \Cref{lem:ExpRV} to the exponential moment of $E(T)$.
	Thus, this restriction comes from the moment estimate used in the proof.
	As our theoretical result indicates, when the stability index $\alpha$ approaches $1$, the convergence rate tends to the standard order of $1/2$, which coincides with the convergence rates established by \cite{jum2014strong, deng2020semi} and others for non-equidistant step size methods. 
	Furthermore, if we consider subordinators with a Laplace exponent that includes a positive linear drift term, as studied in \cite{deng2020semi}, the corresponding inverse subordinator exhibits almost sure Lipschitz continuity. 
	Consequently, there exists a deterministic constant $C > 0$ such that the maximum increment satisfies $\Xi_{\Delta t} \le C\Delta t$.
	Substituting this linear bound directly into the error estimates of \Cref{th:discrete_EM_convergence}, naturally leads to the classical strong convergence rate of $1/2$.
	We will verify this observation numerically in \Cref{sec:num}.
\end{remark}

\section{Strong convergence of truncated EM method}\label{sec:TEM}

The global Lipschitz assumption in \Cref{sec:EM} simplifies the analysis but is often too restrictive for practical applications.
In this section, we extend our analysis to the time-changed SDEs \Cref{eq:SDEs} under relaxed conditions on the coefficients.

\begin{assumption}\label{assum:polynomial_growth}
	There exist positive constants $\gamma$ and $L$ such that for all $x, y \in \mathbb{R}^d$,
	\begin{equation*}
		|f(x) - f(y)| \vee |g(x) - g(y)|
		\le L \bigl(1 + |x|^{\gamma} + |y|^{\gamma}\bigr) |x - y|.
	\end{equation*}
\end{assumption}
Note that \Cref{assum:polynomial_growth} implies the polynomial growth condition: for all $x \in \mathbb{R}^d$,
\begin{equation}\label{eq:super_linear_growth}
	|f(x)| \vee |g(x)| \le M \bigl(1 + |x|^{\gamma+1}\bigr),
\end{equation}
where $M$ is a constant depending on $L$, $f(0)$, and $g(0)$.

\begin{assumption}\label{assum:one_side_Lipschitz}
	There exist constants $p > 2$ and $K > 0$ such that for all $x, y \in \mathbb{R}^d$,
	\begin{equation*}
		(x - y)^\top \bigl(f(x) - f(y)\bigr) + \frac{3p-1}{2} |g(x) - g(y)|^2
		\le K |x - y|^2.
	\end{equation*}
\end{assumption}

\begin{assumption}\label{assum:Khasminskii}
	There exist constants $q > 2$ and $K_1 > 0$ such that for all $x \in \mathbb{R}^d$,
	\begin{equation*}
		x^\top f(x) + \frac{3q-1}{2} |g(x)|^2 \le K_1 \bigl(1 + |x |^2\bigr).
	\end{equation*}
\end{assumption}

Under these assumptions, the time-changed SDEs \Cref{eq:SDEs} admit a unique strong solution \cite{mao1991stability, kobayashi2011stochastic, li2025truncated}.
As discussed in \Cref{sec:intro}, since the coefficients $f$ and $g$ may exhibit super-linear growth, the standard EM method is known to be divergent in general \cite{hutzenthaler2011divergence}. 
To overcome this difficulty, we employ the truncated EM method \cite{mao2015truncate,mao2016convergence}.
A related partially truncated EM method for classical SDEs was studied in \cite{yang2022convergence}.

Let $\mu: \mathbb{R}_+ \to \mathbb{R}_+$ be a strictly increasing continuous function such that $\mu(u) \to \infty$ as $u \to \infty$ and
\begin{equation*}
	\sup_{|x| \le u} \bigl(|f(x)| \vee |g(x)|\bigr) \le \mu(u), \quad \forall u \ge 1.
\end{equation*}
Let $\mu^{-1}: [\mu(0), \infty) \to \mathbb{R}_+$ denote the inverse function of $\mu$. 
We select a constant $\hat{\kappa} \ge 1 \vee \mu(1)$ and a strictly decreasing function $\kappa: (0, 1] \to [\mu(1), \infty)$ satisfying
\begin{equation}\label{eq:trunc_k}
	\lim_{\Delta \to 0} \kappa(\Delta) = \infty 
	\quad \text{and} \quad 
	\Delta^{\frac{\alpha-\varepsilon}{4}}   \kappa(\Delta) \le \hat{\kappa}, \quad \forall \Delta \in (0, 1],
\end{equation}
where $\varepsilon \in (0, \alpha)$ is a fixed parameter.
For a given time step $\Delta \in (0, 1]$, we define the truncation mapping $\pi_\Delta: \mathbb{R}^d \to  \bigl\{x\in\mathbb{R}^d:\ |x|\le \mu^{-1}(\kappa(\Delta))\bigr\}$ by
\begin{equation*}
	\pi_\Delta(x) = \bigl(|x|\wedge \mu^{-1}(\kappa(\Delta))\bigr) \frac{x}{|x|},
\end{equation*}
with the convention $0/0=0$. 
The associated truncated coefficients are given by
\begin{equation*}
	f_\Delta(x) = f\bigl(\pi_\Delta(x)\bigr) \quad \text{and} \quad g_\Delta(x) = g\bigl(\pi_\Delta(x)\bigr).
\end{equation*}
By construction, for any $x \in \mathbb{R}^d$, the truncated coefficients satisfy the bound
\begin{equation}\label{eq:trunc_fg}
	|f_\Delta(x)|\vee |g_\Delta(x)|\le \mu\bigl(\mu^{-1}(\kappa(\Delta))\bigr) = \kappa(\Delta).
\end{equation}

Based on a uniform time grid $t_n = n\Delta$, the truncated EM scheme generates a discrete-time approximation $X_{t_n}$ via $X_0 = X(0)$ and
\begin{equation*}
	X_{t_{n+1}} = X_{t_n} + f_\Delta(X_{t_n}) (E(t_{n+1}) - E(t_n)) + g_\Delta(X_{t_n}) (B(E(t_{n+1})) - B(E(t_n))).
\end{equation*}
We introduce two continuous-time interpolations for the discrete process. 
The first is the piecewise constant process
\begin{equation}\label{eq:step_X}
	\tilde{X}(t) = \sum_{n=0}^\infty X_{t_n} \mathbf{1}_{[t_n, t_{n+1})}(t).
\end{equation}
The second is the continuous process defined by the integral equation
\begin{equation}\label{eq:cont_X}
	\bar{X}(t) = X(0) + \int_0^t f_\Delta\bigl(\tilde{X}(s)\bigr) \mathrm{d}E(s) + \int_0^t g_\Delta\bigl(\tilde{X}(s)\bigr) \mathrm{d}B\bigl(E(s)\bigr), \quad t \in [0, T].
\end{equation}
Note that $\bar{X}(t_n) = X_{t_n}$.

Before proceeding to the error analysis, we present three essential lemmas concerning the properties of the truncated coefficients and the exact solution \cite{guo2018note, li2025truncated}.
\begin{lemma}\label{lem:trunc_Khasminskii}
	Under \Cref{assum:Khasminskii}, for any $\Delta \in (0, 1]$ and $x \in \mathbb{R}^d$,
	\begin{equation*}
		x^\top f_\Delta(x) + \frac{3q-1}{2} |g_\Delta(x)|^2 \le \tilde{K}_1 (1 + |x|^2),
	\end{equation*}
	where $\tilde{K}_1 = 2K_1 (1 \vee [\mu^{-1}(\kappa(1))]^{-1})$.
\end{lemma}

\begin{lemma}\label{lem:trunc_fg}
	Under \Cref{assum:polynomial_growth}, for any $\Delta \in (0, 1]$ and $x, y \in \mathbb{R}^d$,
	\begin{equation*}
		|f_\Delta(x) - f_\Delta(y)|\vee |g_\Delta(x) - g_\Delta(y)|\le L (1 + |x|^\gamma + |y|^\gamma) |x - y|.
	\end{equation*}
\end{lemma}

\begin{lemma}\label{lem:true_solution_bound}
	Suppose \Cref{assum:polynomial_growth,assum:Khasminskii} hold. 
	Then, for any $2 \le p < q$, the exact solution satisfies
	\begin{equation*}
		\mathbb{E}\Bigl[\sup_{0 \le t \le T} |X(t)|^p\Bigr] < \infty.
	\end{equation*}
\end{lemma}

The following lemma provides an error estimate between the continuous interpolation $\bar{X}(t)$ and the step process $\tilde{X}(t)$, expressed in terms of the conditional expectation with respect to the Brownian motion.
\begin{lemma}\label{lem:X_Xbar_diff} 
	Let $\Delta \in (0, 1]$ and $\bar{p} \ge 2$. 
	For any $t \in [t_n, t_{n+1})$, we have almost surely
	\begin{equation*}
		\mathbb{E}_B\Bigl[|\bar{X}(t) - \tilde{X}(t)|^{\bar{p}}\Bigr] \le C_{\bar{p}}   \kappa(\Delta)^{\bar{p}} \Bigl(\Xi_{\Delta}^{\bar{p}} + \Xi_{\Delta}^{\bar{p}/2}\Bigr),
	\end{equation*}
	where $\Xi_\Delta$ is defined in \Cref{lem:sup_delta_E}, and $C_{\bar{p}} > 0$ depends only on $\bar{p}$. 
\end{lemma}

\begin{proof}
	By the definitions of $\bar{X}(t)$ and $\tilde{X}(t)$ in \Cref{eq:cont_X} and \Cref{eq:step_X}, for $t \in [t_n, t_{n+1})$, we have
	\begin{align}\label{eq:X_decomp}
		\mathbb{E}_B\bigl[|\bar{X}(t)-\tilde{X}(t)|^{\bar{p}}\bigr]
		&=\mathbb{E}_B\bigl[|\bar{X}(t)-X_{t_n}|^{\bar{p}}\bigr] \nonumber \\
		&\le 2^{\bar{p}-1}\mathbb{E}_B\biggl[\Bigl|\int_{t_n}^t f_\Delta\bigl(\tilde{X}(s)\bigr)\mathrm{d}E(s)\Bigr|^{\bar{p}}\biggr]
		+2^{\bar{p}-1}\mathbb{E}_B\biggl[\Bigl|\int_{t_n}^t g_\Delta\bigl(\tilde{X}(s)\bigr)\mathrm{d}B\bigl(E(s)\bigr)\Bigr|^{\bar{p}}\biggr].
	\end{align}
	For the first term on the right-hand side, using H\"older's inequality, we have
	\begin{equation}\label{eq:f_Delta_holder}
		\mathbb{E}_B\biggl[\Bigl|\int_{t_n}^t f_\Delta\bigl(\tilde{X}(s)\bigr)\mathrm{d}E(s)\Bigr|^{\bar{p}}\biggr]
		\le (\Delta E_n)^{\bar{p}-1}\mathbb{E}_B\biggl[\int_{t_n}^t\bigl|f_\Delta\bigl(\tilde{X}(s)\bigr)\bigr|^{\bar{p}}\mathrm{d}E(s)\biggr].
	\end{equation}
	For the second term, applying the BDG inequality yields
	\begin{equation}\label{eq:g_Delta_BDG}
		\mathbb{E}_B\biggl[\Bigl|\int_{t_n}^t g_\Delta\bigl(\tilde{X}(s)\bigr)\mathrm{d}B\bigl(E(s)\bigr)\Bigr|^{\bar{p}}\biggr]\le b_{\bar{p}}(\Delta E_n)^{\frac{\bar{p}}{2}-1}\mathbb{E}_B\biggl[\int_{t_n}^t\bigl|g_\Delta\bigl(\tilde{X}(s)\bigr)\bigr|^{\bar{p}}\mathrm{d}E(s)\biggr].
	\end{equation}
	Substituting the estimates \Cref{eq:f_Delta_holder} and \Cref{eq:g_Delta_BDG} into \Cref{eq:X_decomp} and using \Cref{eq:trunc_fg}, we obtain
	\begin{align*}
		\mathbb{E}_B\bigl[|\bar{X}(t)-\tilde{X}(t)|^{\bar{p}}\bigr]
		&\le C_{\bar{p}}\bigl(\kappa(\Delta)\bigr)^{\bar{p}}\Bigl((\Delta E_n)^{\bar{p}}+(\Delta E_n)^{\frac{\bar{p}}{2}}\Bigr) \\
		&\le C_{\bar{p}}\bigl(\kappa(\Delta)\bigr)^{\bar{p}}\Bigl(\Xi_{\Delta}^{\bar{p}}+\Xi_{\Delta}^{\frac{\bar{p}}{2}}\Bigr),
	\end{align*}
	where $C_{\bar{p}}=2^{\bar{p}-1}(b_{\bar{p}} \vee 1)$. 
	This completes the proof.
\end{proof}

The following lemma establishes the uniform moment boundedness of the continuous-time truncated EM solution.
\begin{lemma}\label{lem:Xbar_moment_bound} 
	Under \Cref{assum:polynomial_growth,assum:Khasminskii}, for any $T>0$ and $p \ge 2$, we have
	\begin{equation*}
		\sup_{0<\Delta\le 1}\mathbb{E}\Bigl[\sup_{0\le t\le T}|\bar{X}(t)|^p\Bigr]\le C,
	\end{equation*}
	where $C$ is independent of $\Delta$.
\end{lemma}
\begin{proof}
	Define the stopping time $\zeta_\ell:=\inf\{t\ge 0;|\bar{X}(t)|>\ell\}$ for any positive integer $\ell$. 
	Fix any $\Delta\in(0,1]$ and $T\ge 0$. 
	By It\^o's formula, we deduce from \Cref{eq:cont_X} that, for any $0\le u\le t\wedge\zeta_\ell$,
	\begin{equation*}
		|\bar{X}(u)|^p=|\bar{X}(0)|^p+A_u+M_u,
	\end{equation*}
	where
	\begin{align*}
		A_u&:=\int_0^u\Bigl(p|\bar{X}(s)|^{p-2}\bar{X}^\top(s)f_\Delta\bigl(\tilde{X}(s)\bigr)+\frac{1}{2}p(p-1)|\bar{X}(s)|^{p-2}\bigl|g_\Delta\bigl(\tilde{X}(s)\bigr)\bigr|^2\Bigr)\mathrm{d}E(s),\\
		M_u&:=\int_0^u p|\bar{X}(s)|^{p-2}\bar{X}^\top(s)g_\Delta\bigl(\tilde{X}(s)\bigr)\mathrm{d}B\bigl(E(s)\bigr).
	\end{align*}
	Observe that the stochastic integral $(M_t)_{t\ge 0}$ is a local martingale. 
	By evaluating its quadratic variation and applying Young's inequality, we obtain that for $u\le t\wedge\zeta_\ell$,
	\begin{align*}
		([M,M]_u)^{\frac{1}{2}}
		&\le p\Bigl(\sup_{0\le s\le t\wedge\zeta_\ell}|\bar{X}(s)|^p\int_0^u|\bar{X}(s)|^{p-2}\bigl|g_\Delta\bigl(\tilde{X}(s)\bigr)\bigr|^2\mathrm{d}E(s)\Bigr)^{\frac{1}{2}}\\
		&\le p\Bigl(\frac{\sup_{0\le s\le t\wedge\zeta_\ell}|\bar{X}(s)|^p}{2pb_1}+2pb_1\int_0^u|\bar{X}(s)|^{p-2}\bigl|g_\Delta\bigl(\tilde{X}(s)\bigr)\bigr|^2\mathrm{d}E(s)\Bigr),
	\end{align*}
	where $b_1$ is the constant appearing in the BDG inequality.
	Therefore, by \Cref{lem:trunc_Khasminskii} and Young's inequality,	we obtain
	\begin{align*}
		&\quad \mathbb{E}_B\Bigl[\sup_{0\le u\le t\wedge\zeta_\ell}|\bar{X}(u)|^p\Bigr] \\
		&\le |\bar{X}(0)|^p+\mathbb{E}_B\biggl[\sup_{0\le u\le t\wedge\zeta_\ell}\int_0^u p|\bar{X}(s)|^{p-2}\Bigl(\bar{X}^\top(s)f_\Delta\bigl(\tilde{X}(s)\bigr) +  \frac{1}{2}(p-1)\bigl|g_\Delta\bigl(\tilde{X}(s)\bigr)\bigr|^2\Bigr)\mathrm{d}E(s)\biggr] \\
		&\quad +b_1\mathbb{E}_B\biggl[\frac{1}{2b_1}\sup_{0\le u\le t\wedge\zeta_\ell}|\bar{X}(u)|^p+\int_0^{t\wedge\zeta_\ell} 2b_1p^2|\bar{X}(s)|^{p-2}\bigl|g_\Delta\bigl(\tilde{X}(s)\bigr)\bigr|^2 \mathrm{d}E(s)\biggr]\\
		&\le |\bar{X}(0)|^p+\frac{1}{2}\mathbb{E}_B\Bigl[\sup_{0\le u\le t\wedge\zeta_\ell}|\bar{X}(u)|^p\Bigr] \\
		&\quad +b\mathbb{E}_B\biggl[\sup_{0\le u\le t\wedge\zeta_\ell}\int_0^u p|\bar{X}(s)|^{p-2}\Bigl(\tilde{X}^\top(s)f_\Delta\bigl(\tilde{X}(s)\bigr)+\frac{3p-1}{2}\bigl|g_\Delta\bigl(\tilde{X}(s)\bigr)\bigr|^2\Bigr)\mathrm{d}E(s)\biggr] \\
		&\quad +\mathbb{E}_B\biggl[\sup_{0\le u\le t\wedge\zeta_\ell}\int_0^{u} p|\bar{X}(s)|^{p-2}\bigl(\bar{X}(s)-\tilde{X}(s)\bigr)^\top f_\Delta\bigl(\tilde{X}(s)\bigr)\mathrm{d}E(s)\biggr]\\
		&\le |\bar{X}(0)|^p+\frac{1}{2}\mathbb{E}_B\Bigl[\sup_{0\le u\le t\wedge\zeta_\ell}|\bar{X}(u)|^p\Bigr]+bp\tilde{K}_1\mathbb{E}_B\biggl[\sup_{0\le u\le t\wedge\zeta_\ell}\int_0^u|\bar{X}(s)|^{p-2}(1+|\tilde{X}(s)|^2)\mathrm{d}E(s)\biggr] \\
		&\quad +\mathbb{E}_B\biggl[\sup_{0\le u\le t\wedge\zeta_\ell}\Bigl((p-2)\int_0^u|\bar{X}(s)|^p \mathrm{d}E(s)+2\int_0^u|\bar{X}(s)-\tilde{X}(s)|^{\frac{p}{2}}\bigl|f_\Delta\bigl(\tilde{X}(s)\bigr)\bigr|^{\frac{p}{2}}\mathrm{d}E(s)\Bigr)\biggr],
	\end{align*}
	where $b = 2b_1^2 \vee 1$.
	Consequently, we have
	\begin{align*}
		\mathbb{E}_B\Bigl[\sup_{0\le s\le t\wedge\zeta_\ell}|\bar{X}(s)|^p\Bigr]
		&\le 2|\bar{X}(0)|^p+2bp\tilde{K}_1\mathbb{E}_B\Bigl[\int_0^t|\bar{X}(s\wedge\zeta_\ell)|^{p-2}\bigl(1+|\tilde{X}(s\wedge\zeta_\ell)|^2\bigr)\mathrm{d}E(s)\Bigr]\\
		&\quad+2(p-2)\int_0^t\mathbb{E}_B\bigl[|\bar{X}(s\wedge\zeta_\ell)|^p\bigr]\mathrm{d}E(s)\\
		&\quad+4\int_0^t\mathbb{E}_B\Bigl[|\bar{X}(s\wedge\zeta_\ell)-\tilde{X}(s\wedge\zeta_\ell)|^{\frac{p}{2}}\bigl|f_\Delta\bigl(\tilde{X}(s\wedge\zeta_\ell)\bigr)\bigr|^{\frac{p}{2}}\Bigr]\mathrm{d}E(s).
	\end{align*}
	By \Cref{eq:trunc_k}, \Cref{eq:trunc_fg}, and \Cref{lem:X_Xbar_diff}, it follows that
	\begin{align*}
		& \quad \int_0^t\mathbb{E}_B\Bigl[|\bar{X}(s\wedge\zeta_\ell)-\tilde{X}(s\wedge\zeta_\ell)|^{\frac{p}{2}}\bigl|f_\Delta\bigl(\tilde{X}(s\wedge\zeta_\ell)\bigr)\bigr|^{\frac{p}{2}}\Bigr]\mathrm{d}E(s) \\
		&\le \bigl(\kappa(\Delta)\bigr)^{\frac{p}{2}}\int_0^t \mathbb{E}_B\bigl[|\bar{X}(s\wedge\zeta_\ell)-\tilde{X}(s\wedge\zeta_\ell)|^{\frac{p}{2}}\bigr]\mathrm{d}E(s) \\
		&\le \bigl(\kappa(\Delta)\bigr)^{\frac{p}{2}}\int_0^t \bigl(\mathbb{E}_B\bigl[|\bar{X}(s\wedge\zeta_\ell)-\tilde{X}(s\wedge\zeta_\ell)|^p\bigr]\bigr)^{\frac{1}{2}}\mathrm{d}E(s) \\
		&\le C^{\frac{1}{2}}_{p}\bigl(\kappa(\Delta)\bigr)^{p}E(t)\bigl(\Xi_{\Delta}^{\frac{p}{2}}+\Xi_{\Delta}^{\frac{p}{4}}\bigr).
	\end{align*}
	
	We therefore conclude that
	\begin{align*}
		\mathbb{E}_B\Bigl[\sup_{0\le u\le t\wedge\zeta_\ell}|\bar{X}(u)|^p\Bigr]
		&\le C_1+C_2\int_0^t\Bigl(\mathbb{E}_B\bigl[|\bar{X}(s\wedge\zeta_\ell)|^p\bigr]+\mathbb{E}_B\bigl[|\tilde{X}(s\wedge\zeta_\ell)|^p\bigr]\Bigr)\mathrm{d}E(s) \\
		&\le C_1+2C_2\int_0^t\mathbb{E}_B\Bigl[\sup_{0\le u\le s\wedge\zeta_\ell}|\bar{X}(u)|^p\Bigr]\mathrm{d}E(s),
	\end{align*}
	where
	\begin{equation*}
		C_1=2|\bar{X}(0)|^p+4C^{\frac{1}{2}}_{p}\bigl(\kappa(\Delta)\bigr)^{p}E(T)\bigl(\Xi_{\Delta}^{\frac{p}{2}}+\Xi_{\Delta}^{\frac{p}{4}}\bigr),\quad C_2=2bp\tilde{K}_1\vee 2(p-2).
	\end{equation*}
	Applying Gronwall's inequality, for any $t\in[0,T]$,
	\begin{equation*}
		\mathbb{E}_B\Bigl[\sup_{0\le u\le t\wedge\zeta_\ell}|\bar{X}(u)|^p\Bigr]\le \Bigl(2|\bar{X}(0)|^p+4c^{\frac{1}{2}}_{2p}\bigl(\kappa(\Delta)\bigr)^{p}E(T)\bigl(\Xi_{\Delta}^{\frac{p}{2}}+\Xi_{\Delta}^{\frac{p}{4}}\bigr)\Bigr)e^{2C_2E(t)}.
	\end{equation*}
	Noting that $\zeta_\ell\to\infty$ as $\ell\to\infty$, we set $t=T$ and let $\ell\to\infty$. 
	Taking the expectation $\mathbb{E}_D$ on both sides and applying the Cauchy--Schwarz inequality combined with \Cref{lem:E_moment,lem:ExpRV,lem:sup_delta_E}, we obtain
	\begin{equation*}
		\mathbb{E}\Bigl[\sup_{0\le t\le T}|\bar{X}(t)|^p\Bigr]\le C.
	\end{equation*}
	As this holds for any $\Delta\in(0,1]$ and $C$ is independent of $\Delta$, we have
	\begin{equation*}
		\sup_{0<\Delta\le 1}\mathbb{E}\Bigl[\sup_{0\le t\le T}|\bar{X}(t)|^p\Bigr]\le C.
	\end{equation*}
	The proof is complete. 
\end{proof}

Having established the necessary a priori bounds and local error estimates, we now state the main theorem concerning the strong convergence of the truncated EM scheme.
\begin{theorem}\label{thm:TEM_convergence}
	Let $\bar{p}\in[2,p)$ and $\varepsilon\in(0,\alpha)$.
	Suppose \Cref{assum:polynomial_growth,assum:one_side_Lipschitz,assum:Khasminskii} hold with an exponent $q$ satisfying
	\begin{equation*}
		q>(\gamma+1)p, \qquad \frac{\varepsilon\bigl(q-(\gamma+1)\bar{p}\bigr)} {4(\gamma+2)} \ge \frac{\bar{p}}{2}(\alpha-\varepsilon).
	\end{equation*}
	Then for any $\Delta\in(0,1]$,
	\begin{align*}
		\mathbb{E}\Bigl[\sup_{0\le t\le T}|X(t)-\bar{X}(t)|^{\bar{p}}\Bigr] &\le C\Delta^{\frac{\bar{p}}{2}(\alpha-\varepsilon)},\\
		\mathbb{E}\Bigl[\sup_{0\le t\le T}|X(t)-\tilde{X}(t)|^{\bar{p}}\Bigr]&\le C\Delta^{\frac{\bar{p}}{2}(\alpha-\varepsilon)}.
	\end{align*}
\end{theorem}

\begin{proof}
	Fix $\Delta\in(0,1]$ arbitrarily. 
	Let $e_\Delta(u)=X(u)-\bar{X}(u)$ for $u\ge 0$. 
	For each integer $\ell>|X(0)|$, define the stopping time
	\begin{equation*}
		\theta_{\Delta,\ell}=\inf\{t\ge 0:|X(t)|\vee|\bar{X}(t)|\ge\ell\},
	\end{equation*}
	where we set $\inf\emptyset=\infty$ (as usual, $\emptyset$ denotes the empty set). 
	The dependence of $\theta_{\Delta,\ell}$ on $\Delta$ arises from the numerical interpolation $\bar{X}$.
	Following a similar procedure to the proof of \Cref{lem:Xbar_moment_bound}, by applying the It\^{o} formula to $|e_\Delta(u)|^{\bar{p}}$, we can obtain
	\begin{align}\label{eq:error_BDG}
		&\quad \mathbb{E}_B\Bigl[\sup_{0\le u \le t\wedge\theta_{\Delta,\ell}}|e_\Delta(u)|^{\bar{p}}\Bigr] \nonumber \\
		&\le b\mathbb{E}_B\Bigl[\sup_{0\le u\le t\wedge\theta_{\Delta,\ell}}\int_0^{u}\bar{p}|e_\Delta(s)|^{\bar{p}-2}\Bigl(e_\Delta^\top(s)\bigl(f\bigl(X(s)\bigr)-f_\Delta\bigl(\tilde{X}(s)\bigr)\bigr)+\frac{3\bar{p}-1}{2}\bigl|g\bigl(X(s)\bigr)-g_\Delta\bigl(\tilde{X}(s)\bigr)\bigr|^2\Bigr)\mathrm{d}E(s)\Bigr] \nonumber \\
		&\quad+\frac{1}{2}\mathbb{E}_B\Bigl[\sup_{0\le u\le t\wedge\theta_{\Delta,\ell}}|e_\Delta(u)|^{\bar{p}}\Bigr].
	\end{align}
	Applying Young's inequality, we can deduce that
	\begin{align*}
		&\quad\frac{3\bar{p}-1}{2}\bigl|g\bigl(X(s)\bigr)-g_\Delta\bigl(\tilde{X}(s)\bigr)\bigr|^2\\
		&\le\frac{3\bar{p}-1}{2}\Bigl(\bigl|g\bigl(X(s)\bigr)-g\bigl(\bar{X}(s)\bigr)\bigr|^2+2\bigl|g\bigl(X(s)\bigr)-g\bigl(\bar{X}(s)\bigr)\bigr|\bigl|g\bigl(\bar{X}(s)\bigr)-g_\Delta\bigl(\tilde{X}(s)\bigr)\bigr|+\bigl|g\bigl(\bar{X}(s)\bigr)-g_\Delta\bigl(\tilde{X}(s)\bigr)\bigr|^2\Bigr)\\
		&\le\frac{3\bar{p}-1}{2}\Bigl[\bigl|g\bigl(X(s)\bigr)-g\bigl(\bar{X}(s)\bigr)\bigr|^2+\Bigl(\frac{3p-3\bar{p}}{3\bar{p}-1}\bigl|g\bigl(X(s)\bigr)-g\bigl(\bar{X}(s)\bigr)\bigr|^2+\frac{3\bar{p}-1}{3p-3\bar{p}}\bigl|g\bigl(\bar{X}(s)\bigr)-g_\Delta\bigl(\tilde{X}(s)\bigr)\bigr|^2\Bigr)\\
		&\quad+\bigl|g\bigl(\bar{X}(s)\bigr)-g_\Delta\bigl(\tilde{X}(s)\bigr)\bigr|^2\Bigr]\\
		&=\frac{3p-1}{2}\bigl|g\bigl(X(s)\bigr)-g\bigl(\bar{X}(s)\bigr)\bigr|^2+\frac{(3\bar{p}-1)(3p-1)}{2(3p-3\bar{p})}\bigl|g\bigl(\bar{X}(s)\bigr)-g_\Delta\bigl(\tilde{X}(s)\bigr)\bigr|^2.
	\end{align*}
	Substituting this into \Cref{eq:error_BDG}, we have
	\begin{equation}\label{eq:error_decomp}
		\mathbb{E}_B\Bigl[\sup_{0\le u\le t\wedge\theta_{\Delta,\ell}}|e_\Delta(u)|^{\bar{p}}\Bigr]
		\le 2\mathbb{E}_B\Bigl[\sup_{0\le u\le t\wedge\theta_{\Delta,\ell}}J_1(u)\Bigr] + 2\mathbb{E}_B\Bigl[\sup_{0\le u\le t\wedge\theta_{\Delta,\ell}}J_2(u)\Bigr],
	\end{equation}
	where
	\begin{align*}
		J_1(u)&:=b\int_0^{u}\bar{p}|e_\Delta(s)|^{\bar{p}-2}\Bigl(e_\Delta^\top(s)\bigl(f(X(s))-f(\bar{X}(s))\bigr)+\frac{3p-1}{2}\bigl|g\bigl(X(s)\bigr)-g\bigl(\bar{X}(s)\bigr)\bigr|^2\Bigr)\mathrm{d}E(s),\\
		J_2(u)&:=b\int_0^{u}\bar{p}|e_\Delta(s)|^{\bar{p}-2}\Bigl(e_\Delta^\top(s)\bigl(f(\bar{X}(s))-f_\Delta(\tilde{X}(s))\bigr)+\frac{(3\bar{p}-1)(3p-1)}{2(3p-3\bar{p})}\bigl|g(\bar{X}(s))-g_\Delta(\tilde{X}(s))\bigr|^2\Bigr)\mathrm{d}E(s).
	\end{align*}
	
	For $J_1(u)$, by \Cref{assum:one_side_Lipschitz}, we have
	\begin{equation}\label{eq:J1_bound}
		J_1(u)\le C_1\int_0^{u}|e_\Delta(s)|^{\bar{p}}\mathrm{d}E(s),
	\end{equation}
	where $C_1=b\bar{p}K$. 
	
	We can decompose $J_2(u)$ into two parts as follows
	\begin{align}\label{eq:J2_split}
		J_2(u)&\le b\int_0^{u}\bar{p}|e_\Delta(s)|^{\bar{p}-2}\Bigl(e_\Delta^\top(s)\bigl(f(\bar{X}(s))-f(\tilde{X}(s))\bigr)+\frac{(3\bar{p}-1)(3p-1)}{3p-3\bar{p}}\bigl|g(\bar{X}(s))
		-g(\tilde{X}(s))\bigr|^2\Bigr)\mathrm{d}E(s)\nonumber\\
		&\quad+b\int_0^{u}\bar{p}|e_\Delta(s)|^{\bar{p}-2}\Bigl(e_\Delta^\top(s)\bigl(f(\tilde{X}(s))-f_\Delta(\tilde{X}(s))\bigr)+\frac{(3\bar{p}-1)(3p-1)}{3p-3\bar{p}}\bigl|g(\tilde{X}(s))
		-g_\Delta(\tilde{X}(s))\bigr|^2\Bigr)\mathrm{d}E(s)\nonumber\\
		&=:J_{21}+J_{22}.
	\end{align}
	For $J_{21}$, by Young's inequality and \Cref{assum:polynomial_growth}, we can show that
	\begin{align}\label{eq:J21_bound}
		&\quad J_{21}\nonumber\\
		&\le b\int_0^{u}\bar{p}|e_\Delta(s)|^{\bar{p}-2}\Bigl(\frac{1}{2}|e_\Delta(s)|^2+\frac{1}{2}\bigl|f\bigl(\bar{X}(s)\bigr)-f\bigl(\tilde{X}(s)\bigr)\bigr|^2 +\frac{(3\bar{p}-1)(3p-1)}{3p-3\bar{p}}\bigl|g\bigl(\bar{X}(s)\bigr)-g\bigl(\tilde{X}(s)\bigr)\bigr|^2\Bigr)\mathrm{d}E(s) \nonumber \\
		&\le \frac{(\bar{p}-1)(3p-3\bar{p})+(\bar{p}-2)(3\bar{p}-1)(3p-1)}{3p-3\bar{p}}b\int_0^{u}|e_\Delta(s)|^{\bar{p}}\mathrm{d}E(s)  +b\int_0^{u}\bigl|f\bigl(\bar{X}(s)\bigr)-f\bigl(\tilde{X}(s)\bigr)\bigr|^{\bar{p}}\mathrm{d}E(s) \nonumber \\
		&\quad +\frac{2(3\bar{p}-1)(3p-1)}{3p-3\bar{p}}b\int_0^{u}\bigl|g\bigl(\bar{X}(s)\bigr)-g\bigl(\tilde{X}(s)\bigr)\bigr|^{\bar{p}}\mathrm{d}E(s) \nonumber \\
		&\le C_{21}\Bigl(\int_0^{u}|e_\Delta(s)|^{\bar{p}}\mathrm{d}E(s)+2L\int_0^{u}\bigl(1+|\bar{X}(s)|^{\gamma\bar{p}}+|\tilde{X}(s)|^{\gamma\bar{p}}\bigr)|\bar{X}(s)-\tilde{X}(s)|^{\bar{p}}\mathrm{d}E(s)\Bigr),
	\end{align}
	where
	\begin{equation*}
		C_{21}=\max\Bigl\{\frac{(\bar{p}-1)(3p-3\bar{p})+(\bar{p}-2)(3\bar{p}-1)(3p-1)}{3p-3\bar{p}},1,\frac{2(3\bar{p}-1)(3p-1)}{3p-3\bar{p}}\Bigr\}b.
	\end{equation*}
	Following a similar derivation for $J_{22}$, we obtain
	\begin{equation}\label{eq:J22_bound}
		J_{22}\le C_{22}\Bigl(\int_0^{u}|e_\Delta(s)|^{\bar{p}}\mathrm{d}E(s)+2L\int_0^{u}\bigl(1+|\tilde{X}(s)|^{\gamma\bar{p}}+\bigl|\pi_\Delta\bigl(\tilde{X}(s)\bigr)\bigr|^{\gamma\bar{p}}\bigr)\bigl|\tilde{X}(s)-\pi_\Delta\bigl(\tilde{X}(s)\bigr)\bigr|^{\bar{p}}\mathrm{d}E(s)\Bigr),
	\end{equation}
	where $C_{22}$ is a positive constant which only depends on $b$, $p$, and $\bar{p}$. 
	
	Substituting \Cref{eq:J1_bound,eq:J2_split,eq:J21_bound,eq:J22_bound} into \Cref{eq:error_decomp}, we obtain
	\begin{align}\label{eq:combined_bound}
		\mathbb{E}_B\Bigl[\sup_{0\le u\le t\wedge\theta_{\Delta,\ell}}|e_\Delta(u)|^{\bar{p}}\Bigr]
		&\le 2(C_1+C_{21}+C_{22})\mathbb{E}_B\biggl[\sup_{0\le u\le t\wedge\theta_{\Delta,\ell}}\int_0^u|e_\Delta(s)|^{\bar{p}}\mathrm{d}E(s)\biggr]+4LI_{21}+4LI_{22},
	\end{align}
	where
	\begin{align*}
		I_{21}&:=C_{21}\mathbb{E}_B\biggl[\sup_{0\le u\le t\wedge\theta_{\Delta,\ell}}\int_0^u\bigl(1+|\bar{X}(s)|^{\gamma\bar{p}}+|\tilde{X}(s)|^{\gamma\bar{p}}\bigr)|\bar{X}(s)-\tilde{X}(s)|^{\bar{p}}\mathrm{d}E(s)\biggr], \\
		I_{22}&:=C_{22}\mathbb{E}_B\biggl[\sup_{0\le u\le t\wedge\theta_{\Delta,\ell}}\int_0^u\bigl(1+|\tilde{X}(s)|^{\gamma\bar{p}}+\bigl|\pi_\Delta\bigl(\tilde{X}(s)\bigr)\bigr|^{\gamma{\bar{p}}}\bigr)\bigl|\tilde{X}(s)-\pi_\Delta\bigl(\tilde{X}(s)\bigr)\bigr|^{\bar{p}}\mathrm{d}E(s)\biggr].
	\end{align*}
	For $I_{21}$, using H\"older's inequality and \Cref{lem:X_Xbar_diff}, and noting that for $r = \bar{p}/q \in (0,1)$ the inequality $(x+y)^r \le x^r + y^r$ holds for any $x, y \ge 0$, we have
	\begin{align}\label{eq:I21_estimate}
		I_{21}	&\le C_{21}\int_0^t\Bigl(\mathbb{E}_B\bigl[(1+|\bar{X}(s\wedge\theta_{\Delta,\ell})|^{\gamma\bar{p}}+|\tilde{X}(s\wedge\theta_{\Delta,\ell})|^{\gamma\bar{p}})^{\frac{q}{q-\bar{p}}}\bigr]\Bigr)^{\frac{q-\bar{p}}{q}}\Bigl(\mathbb{E}_B\bigl[|\bar{X}(s\wedge\theta_{\Delta,\ell})-\tilde{X}(s\wedge\theta_{\Delta,\ell})|^{q}\bigr]\Bigr)^{\frac{\bar{p}}{q}}\mathrm{d}E(s) \nonumber \\
		&\le C_{21}(\kappa(\Delta))^{\bar{p}}\Bigl(\Xi_{\Delta}^{q}+\Xi_{\Delta}^{\frac{q}{2}}\Bigr)^{\frac{\bar{p}}{q}}\int_0^t\Bigl(\mathbb{E}_B\Bigl[(1+|\bar{X}(s\wedge\theta_{\Delta,\ell})|^{\gamma\bar{p}}+|\tilde{X}(s\wedge\theta_{\Delta,\ell})|^{\gamma\bar{p}})^{\frac{q}{q-\bar{p}}}\Bigr]\Bigr)^{\frac{q-\bar{p}}{q}}\mathrm{d}E(s) \nonumber\\ 
		&\le C_{21}\bigl(\kappa(\Delta)\bigr)^{\bar{p}}\bigl(\Xi_{\Delta}^{\bar{p}}+\Xi_{\Delta}^{\frac{\bar{p}}{2}}\bigr)\int_0^t\Bigl(\mathbb{E}_B\bigl[(1+|\bar{X}(s\wedge\theta_{\Delta,\ell})|^{\gamma\bar{p}}+|\tilde{X}(s\wedge\theta_{\Delta,\ell})|^{\gamma\bar{p}})^{\frac{q}{q-\bar{p}}}\bigr]\Bigr)^{\frac{q-\bar{p}}{q}}\mathrm{d}E(s).
	\end{align}
	
	For $I_{22}$, by H\"older's inequality, we obtain
	\begin{align}\label{eq:I22_estimate}
		I_{22}&\le C_{22}\int_0^t\Bigl(\mathbb{E}_B\bigl[1+|\tilde{X}(s\wedge\theta_{\Delta,\ell})|^q+\bigl|\pi_\Delta\bigl(\tilde{X}(s\wedge\theta_{\Delta,\ell})\bigr)\bigr|^{q}\bigr]\Bigr)^{\frac{\gamma\bar{p}}{q}} \Bigl(\mathbb{E}_B\bigl[|\tilde{X}(s\wedge\theta_{\Delta,\ell})-\pi_\Delta\bigl(\tilde{X}(s\wedge\theta_{\Delta,\ell})\bigr)\bigr|^{\frac{q\bar{p}}{q-\gamma\bar{p}}}\bigr]\Bigr)^{\frac{q-\gamma\bar{p}}{q}}\mathrm{d}E(s) \nonumber \\
		&\le C_{22}\int_0^t\Bigl(\mathbb{E}_B\bigl[1+|\tilde{X}(s\wedge\theta_{\Delta,\ell})|^q+\bigl|\pi_\Delta\bigl(\tilde{X}(s\wedge\theta_{\Delta,\ell})\bigr)\bigr|^{q}\bigr]\Bigr)^{\frac{\gamma\bar{p}}{q}}\Bigl(\mathbb{E}_B\bigl[\mathbf{1}_{\{|\tilde{X}(s\wedge\theta_{\Delta,\ell})|>\mu^{-1}(\kappa(\Delta))\}}|\tilde{X}(s\wedge\theta_{\Delta,\ell})|^{\frac{q\bar{p}}{q-\gamma\bar{p}}}\bigr]\Bigr)^{\frac{q-\gamma\bar{p}}{q}}\mathrm{d}E(s) \nonumber \\
		&\le C_{22}\int_0^t\Bigl(\mathbb{E}_B\bigl[1+|\tilde{X}(s\wedge\theta_{\Delta,\ell})|^q+\bigl|\pi_\Delta\bigl(\tilde{X}(s\wedge\theta_{\Delta,\ell})\bigr)\bigr|^{q}\bigr]\Bigr)^{\frac{\gamma\bar{p}}{q}}\nonumber\\
		&\quad \times\Bigl(\mathbb{P}_B\bigl[|\tilde{X}(s\wedge\theta_{\Delta,\ell})|>\mu^{-1}\bigl(\kappa(\Delta)\bigr)\bigr]^{\frac{q-(\gamma+1)\bar{p}}{q-\gamma\bar{p}}}\bigl(\mathbb{E}_B\bigl[|\tilde{X}(s\wedge\theta_{\Delta,\ell})|^q\bigr]\bigr)^{\frac{\bar{p}}{q-\gamma\bar{p}}}\Bigr)^{\frac{q-\gamma\bar{p}}{q}}\mathrm{d}E(s) \nonumber \\
		&\le C_{22}\bigl(\mu^{-1}\bigl(\kappa(\Delta)\bigr)\bigr)^{(\gamma+1)\bar{p}-q}\int_0^t\Bigl(\mathbb{E}_B\bigl[1+|\tilde{X}(s\wedge\theta_{\Delta,\ell})|^q+\bigl|\pi_\Delta\bigl(\tilde{X}(s\wedge\theta_{\Delta,\ell})\bigr)\bigr|^{q}\bigr]\Bigr)^{\frac{\gamma\bar{p}}{q}}\bigl(\mathbb{E}_B\bigl[|\tilde{X}(s\wedge\theta_{\Delta,\ell})|^q\bigr]\bigr)^{\frac{q-\gamma\bar{p}}{q}}\mathrm{d}E(s). 
	\end{align}
	
	Substituting \Cref{eq:I21_estimate} and \Cref{eq:I22_estimate} into \Cref{eq:combined_bound} gives
	\begin{align*}
		\mathbb{E}_B\Bigl[\sup_{0\le u\le t\wedge\theta_{\Delta,\ell}}|e_\Delta(u)|^{\bar{p}}\Bigr]
		&\le 2(C_1+C_{21}+C_{22})\mathbb{E}_B\Bigl[\sup_{0\le u\le t\wedge\theta_{\Delta,\ell}}\int_0^u|e_\Delta(s)|^{\bar{p}}\mathrm{d}E(s)\Bigr] \\
		&\quad+4LC_{21}\bigl(\kappa(\Delta)\bigr)^{\bar{p}}\bigl(\Xi_\Delta^{\frac{\bar{p}}{2}}+\Xi_\Delta^{\bar{p}}\bigr)\Lambda_1 +4LC_{22}\bigl(\mu^{-1}\bigl(\kappa(\Delta)\bigr)\bigr)^{(\gamma+1)\bar{p}-q}\Lambda_2 \\
		&\le 2(C_1+C_{21}+C_{22})\int_0^t\mathbb{E}_B\Bigl[\sup_{0\le r\le s\wedge\theta_{\Delta,\ell}}|e_\Delta(r)|^{\bar{p}}\Bigr]\mathrm{d}E(s) \\
		&\quad+4LC_{21}\bigl(\kappa(\Delta)\bigr)^{\bar{p}}\bigl(\Xi_\Delta^{\frac{\bar{p}}{2}}+\Xi_\Delta^{\bar{p}}\bigr)\Lambda_1 +4LC_{22}\bigl(\mu^{-1}\bigl(\kappa(\Delta)\bigr)\bigr)^{(\gamma+1)\bar{p}-q}\Lambda_2,
	\end{align*}
	where
	\begin{align*}
		\Lambda_1:&=\int_0^T\Bigl(\mathbb{E}_B\bigl[(1+|\bar{X}(s\wedge\theta_{\Delta,\ell})|^{\gamma\bar{p}}+|\tilde{X}(s\wedge\theta_{\Delta,\ell})|^{\gamma\bar{p}})^{\frac{q}{q-\bar{p}}}\bigr]\Bigr)^{\frac{q-\bar{p}}{q}}\mathrm{d}E(s),\\
		\Lambda_2:&=\int_0^T\Bigl(\mathbb{E}_B\bigl[1+|\tilde{X}(s\wedge\theta_{\Delta,\ell})|^q+\bigl|\pi_\Delta(\tilde{X}(s\wedge\theta_{\Delta,\ell}))\bigr|^{q}\bigr]\Bigr)^{\frac{\gamma\bar{p}}{q}}\Bigl(\mathbb{E}_B\bigl[|\tilde{X}(s\wedge\theta_{\Delta,\ell})|^q\bigr]\Bigr)^{\frac{q-\gamma\bar{p}}{q}}\mathrm{d}E(s).
	\end{align*}
	Using Gronwall's inequality yields that
	\begin{equation*}
		\mathbb{E}_B\Bigl[\sup_{0\le u\le t\wedge\theta_{\Delta,\ell}}|e_\Delta(u)|^{\bar{p}}\Bigr]
		\le C\Bigl(\bigl(\kappa(\Delta)\bigr)^{\bar{p}}\bigl(\Xi_\Delta^{\frac{\bar{p}}{2}}+\Xi_\Delta^{\bar{p}}\bigr)\Lambda_1+\bigl(\mu^{-1}\bigl(\kappa(\Delta)\bigr)\bigr)^{(\gamma+1)\bar{p}-q}\Lambda_2\Bigr)e^{C E(T)}.
	\end{equation*}
	Setting $t=T$, taking $\mathbb{E}_D$ on both sides, and applying the Cauchy--Schwarz inequality, we obtain
	\begin{align*}
		\mathbb{E}\Bigl[\sup_{0\le u\le T\wedge\theta_{\Delta,\ell}}|e_\Delta(u)|^{\bar{p}}\Bigr]
		&\le C\bigl(\kappa(\Delta)\bigr)^{\bar{p}}
		\Bigl(\mathbb{E}_D\bigl[\Xi_\Delta^{\bar{p}}+\Xi_\Delta^{2\bar{p}}\bigr]\Bigr)^{\frac12}
		\Bigl(\mathbb{E}_D\bigl[\Lambda_1^2e^{2CE(T)}\bigr]\Bigr)^{\frac12}\nonumber\\
		&\quad+C\bigl(\mu^{-1}\bigl(\kappa(\Delta)\bigr)\bigr)^{(\gamma+1)\bar{p}-q}
		\Bigl(\mathbb{E}_D\bigl[\Lambda_2^2e^{2CE(T)}\bigr]\Bigr)^{\frac12}.
	\end{align*}
	Using \Cref{lem:ExpRV,lem:Xbar_moment_bound}, the relation $|\pi_\Delta(x)|\le |x|$, and the bound $|\tilde{X}(s)|\le\sup_{0\le u\le s}|\bar{X}(u)|$, we can deduce
	\begin{equation*}
		\Bigl(\mathbb{E}_D\bigl[\Lambda_1^2e^{2CE(T)}\bigr]\Bigr)^{\frac12},
		\qquad
		\Bigl(\mathbb{E}_D\bigl[\Lambda_2^2e^{2CE(T)}\bigr]\Bigr)^{\frac12}
	\end{equation*}
	are bounded uniformly in $\Delta$ and $\ell$. 
	For the fixed $\Delta$, it follows from \Cref{lem:true_solution_bound,lem:Xbar_moment_bound} that $\theta_{\Delta,\ell}\rightarrow\infty$ almost surely as $\ell\to\infty$, and hence $T\wedge\theta_{\Delta,\ell}\rightarrow T$ almost surely. 
	Therefore, letting $\ell\to\infty$, applying the monotone convergence theorem, and using \Cref{lem:sup_delta_E} (with $\varepsilon/2$ in place of $\varepsilon$) yield
	\begin{align}\label{eq:error_bound_kuppa}
		\mathbb{E}\Bigl[\sup_{0\le t\le T}|X(t)-\bar{X}(t)|^{\bar{p}}\Bigr]
		&\le C\bigl(\kappa(\Delta)\bigr)^{\bar{p}}\Bigl(\mathbb{E}_D\bigl[\Xi_\Delta^{\bar{p}}+\Xi_\Delta^{2\bar{p}}\bigr]\Bigr)^{\frac12}
		+C\bigl(\mu^{-1}\bigl(\kappa(\Delta)\bigr)\bigr)^{(\gamma+1)\bar{p}-q}\nonumber\\
		&\le C\Bigl(\bigl(\kappa(\Delta)\bigr)^{\bar{p}}\Delta^{\frac{\bar{p}}{2}(\alpha-\frac{\varepsilon}{2})}+\bigl(\mu^{-1}\bigl(\kappa(\Delta)\bigr)\bigr)^{(\gamma+1)\bar{p}-q}\Bigr).
	\end{align}	
	Following the truncation framework, we specify the parameters as
	\begin{equation*}
		\kappa(\Delta)=\Delta^{-\frac{\varepsilon}{4}}, \quad  \mu(u) = 2Mu^{\gamma + 2},
	\end{equation*}
	where $\varepsilon>0$ is chosen to be sufficiently small such that $\kappa(\Delta)$ is consistent with \Cref{eq:trunc_k}.
	Substituting these choices into \Cref{eq:error_bound_kuppa}, we arrive at
	\begin{align*}
		\mathbb{E}\Bigl[\sup_{0\le t\le T}|X(t)-\bar{X}(t)|^{\bar{p}}\Bigr]
		&\le C\Bigl(\Delta^{\frac{\bar{p}}{2}(\alpha - \varepsilon)}+\Delta^{\frac{\varepsilon(q-(\gamma+1)\bar{p})}{4(\gamma+2)}}\Bigr),\\
		\mathbb{E}\Bigl[\sup_{0\le t\le T}|X(t)-\tilde{X}(t)|^{\bar{p}}\Bigr]
		&\le C\Bigl(\Delta^{\frac{\bar{p}}{2}(\alpha - \varepsilon)}+\Delta^{\frac{\varepsilon(q-(\gamma+1)\bar{p})}{4(\gamma+2)}}\Bigr).
	\end{align*}
	By the condition on $q$ in \Cref{thm:TEM_convergence}, we have
	\begin{equation*}
		\frac{\varepsilon(q-(\gamma+1)\bar{p})}{4(\gamma+2)} \ge \frac{\bar{p}}{2}(\alpha - \varepsilon).
	\end{equation*}
	Hence, the second term has order at least $\Delta^{\frac{\bar{p}}{2}(\alpha-\varepsilon)}$,	which completes the proof.
\end{proof}

\section{Numerical simulations}\label{sec:num}

In this section, we present numerical experiments to illustrate our theoretical results. 
We first define our error metric and then apply our numerical scheme to several examples.

To quantify the accuracy of our numerical scheme, we employ the $L_1$ error metric defined as
\begin{equation*}
	e_i = \mathbb{E}\Bigl[\sup_{0\le t\le T}|X(t) - X_{\Delta t_i}(t)|\Bigr],
\end{equation*}
where $X_{\Delta t_i}(t)$ represents the numerical solution at time $t$ computed with step size $\Delta t_i = 2^{-i}$, and $X(t)$ denotes the true solution. 
For our experiments, we set $T=1$ and approximate the expectation via Monte Carlo simulation with $M=10^3$ independent samples
\begin{equation*}
	e_i \approx \frac{1}{M}\sum_{j=1}^{M}\sup_{0\le t\le T}\Bigl|X^j(t) - X^j_{\Delta t_i}(t)\Bigr|.
\end{equation*}
We employ a sufficiently fine discretization with step size $\Delta t = 2^{-20}$ as the reference solution.
The convergence analysis is conducted by computing numerical solutions with a sequence of decreasing step sizes $\{2^{10}\Delta t, 2^9\Delta t, 2^8\Delta t, 2^7\Delta t\}$. 
The corresponding convergence rates are then determined through linear regression in logarithmic coordinates, specifically by fitting $\log(e_i)$ against $\log(\Delta t_i)$.

\begin{example}\label{ex:1}
	We consider a two-dimensional coupled time-changed SDE system
	\begin{equation*}
		\left\{
		\begin{aligned}
			\mathrm{d}X_1(t) &= -(X_1(t) + X_2(t))  \mathrm{d}E(t) + 2(X_1(t) + X_2(t))  \mathrm{d}B_1(E(t)),  \\
			\mathrm{d}X_2(t) &= -2(X_1(t) + X_2(t))  \mathrm{d}E(t) + (X_1(t) + X_2(t))  \mathrm{d}B_2(E(t)), 	
		\end{aligned}
		\right.
	\end{equation*}
	with initial conditions $X_1(0) = 1$ and $X_2(0) = 2$. Here, $B_1$ and $B_2$ are independent standard Brownian motions.
\end{example}

\begin{figure}[htbp!]
	\centering
	\includegraphics[width=0.48\linewidth]{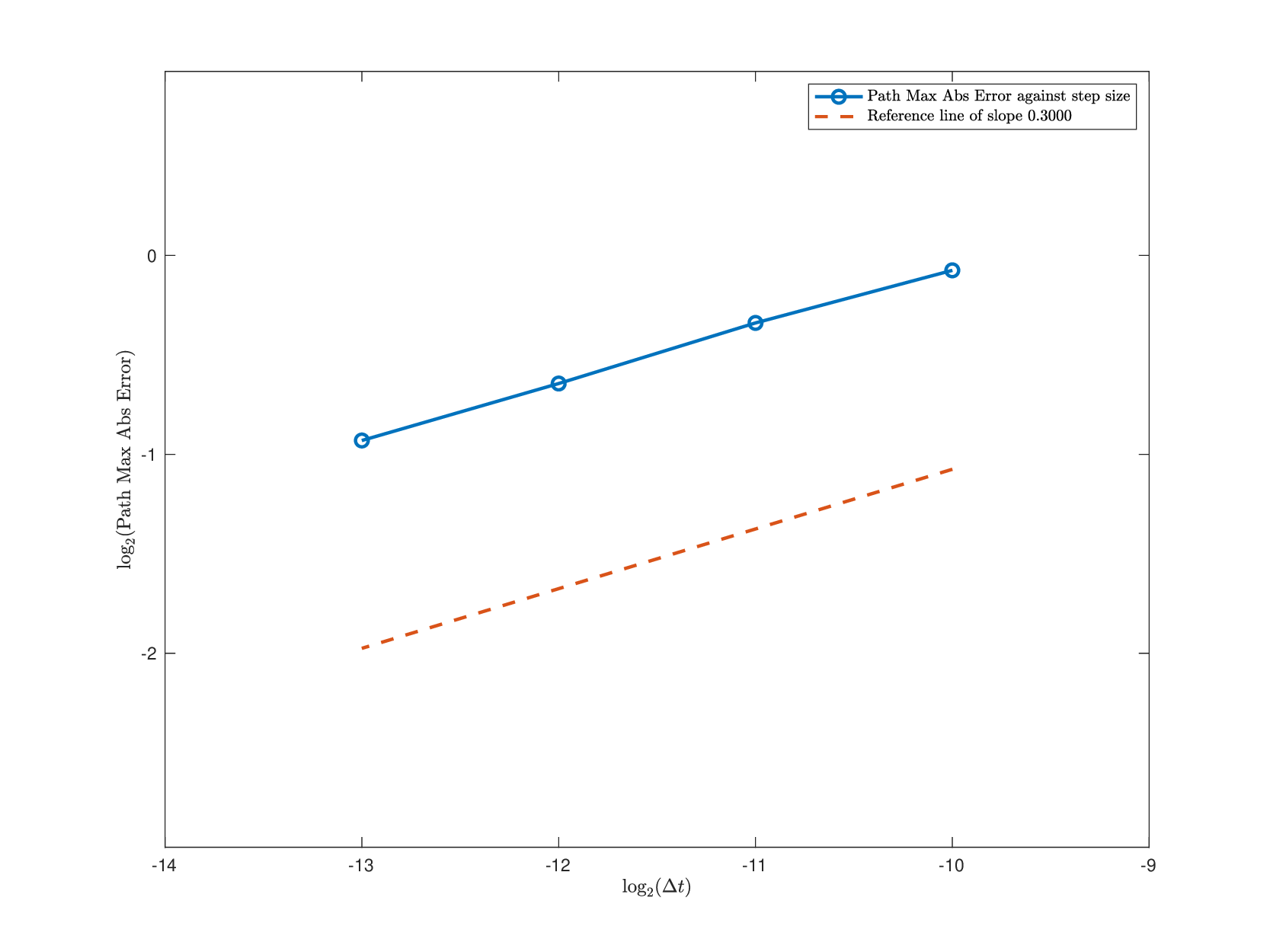}
	\hfill
	\includegraphics[width=0.48\linewidth]{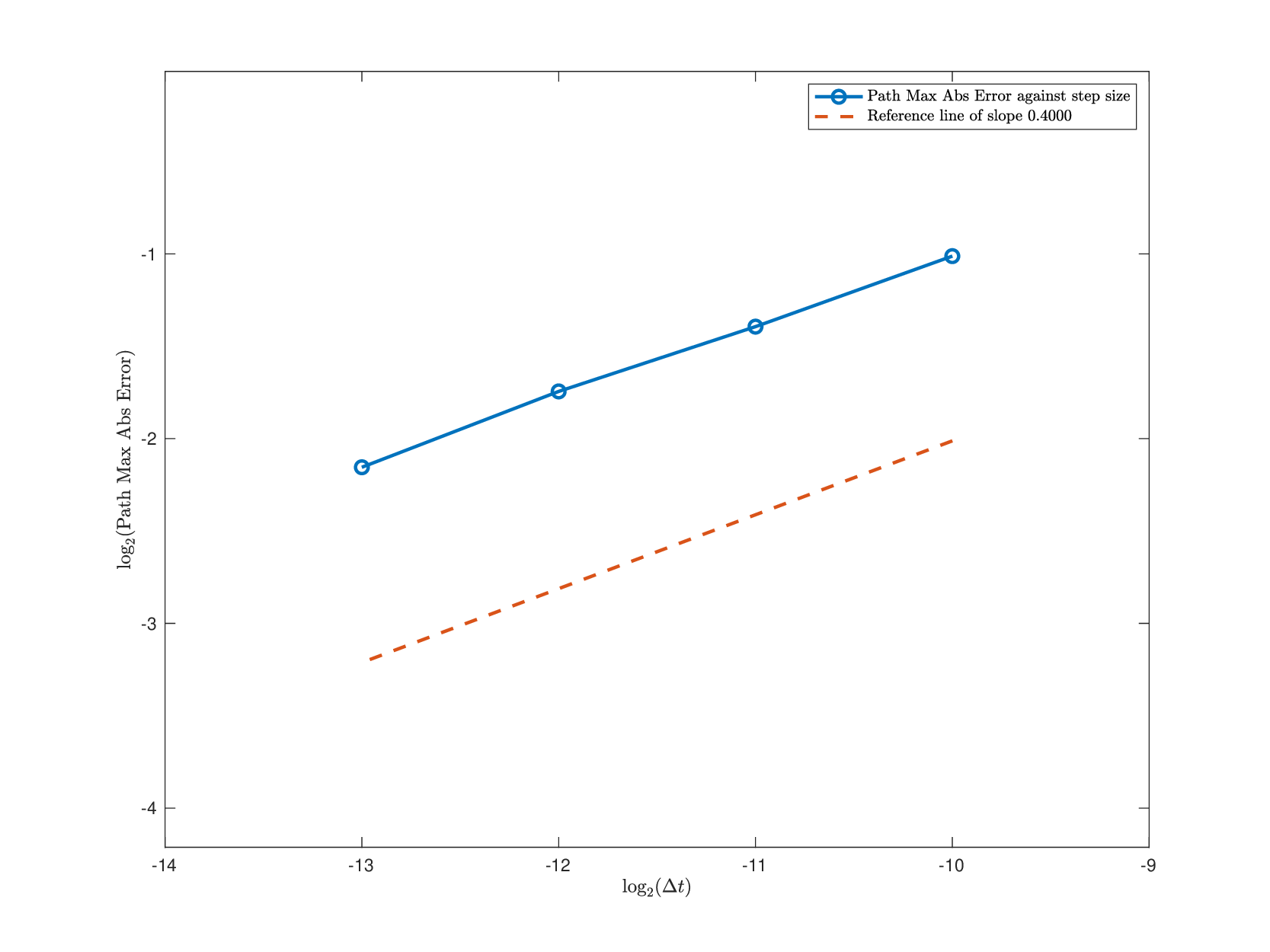}
	\caption{Convergence rates for $\alpha=0.6$ (left) and $\alpha=0.8$ (right) in \Cref{ex:1}.}
	\label{fig:convergence_alpha}
\end{figure}

One can easily verify that this example satisfies all assumptions of \Cref{th:cotinuous_EM_convergence}.
We calculate the numerical convergence rates using the method described above.
The numerical results for this coupled example align well with our theoretical predictions, thus providing strong support for the validity of our analytical framework.
Additionally, we provide a comparison of the convergence rates for $\alpha=0.6$ and $\alpha=0.8$ in \Cref{fig:convergence_alpha}, which visually demonstrates the agreement between theoretical and numerical results.

\begin{figure}[htbp!]
	\centering
	\includegraphics[width=0.48\linewidth]{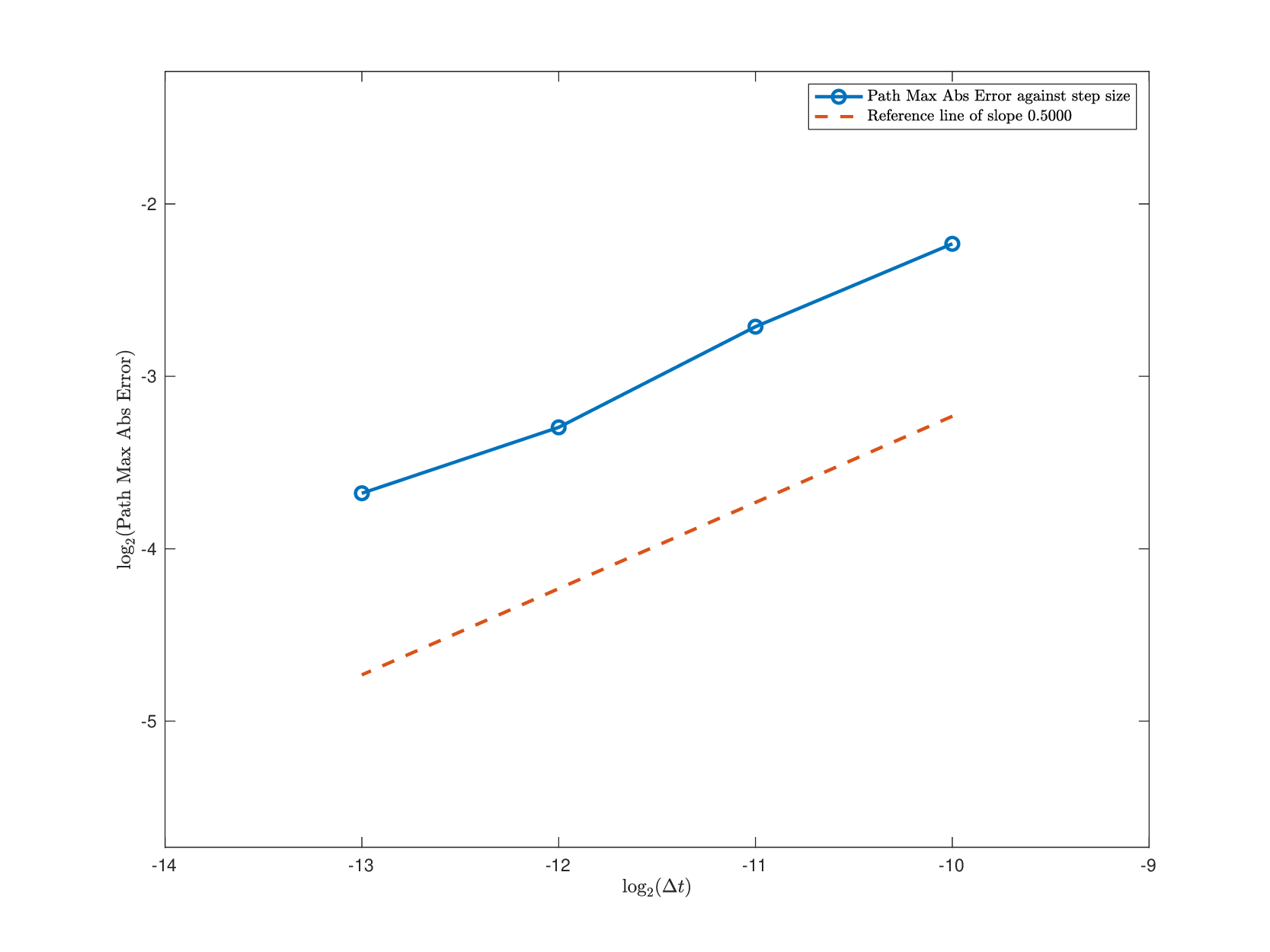}
	\hfill
	\includegraphics[width=0.48\linewidth]{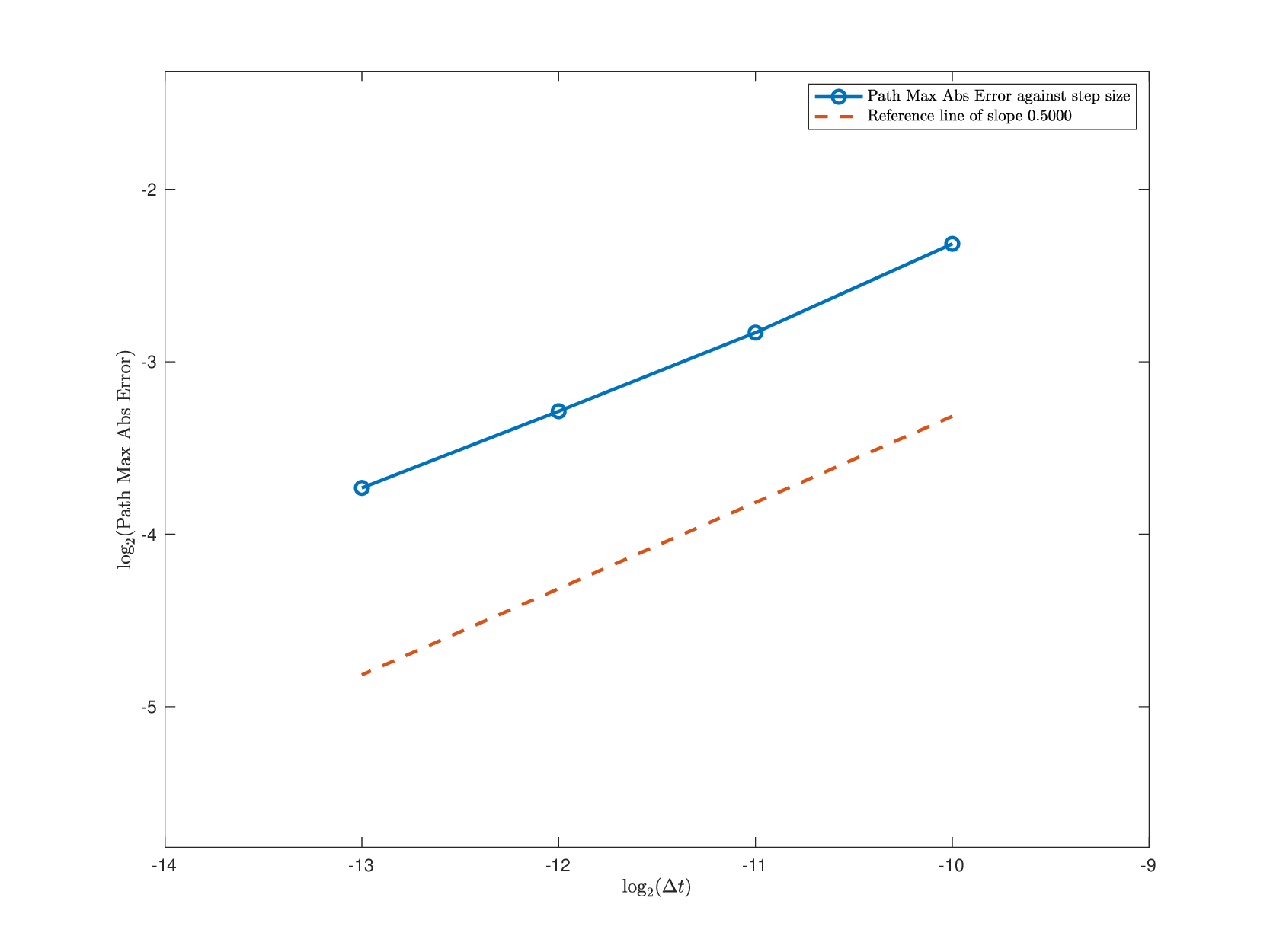}
	\caption{Convergence rates with positive linear drift coefficient $1$ for $\alpha=0.6$ (left) and $\alpha=0.8$ (right) in \Cref{ex:1}.}
	\label{fig:convergence_with_drift}
\end{figure}

To verify the observation in the remark following \Cref{th:cotinuous_EM_convergence}, we add a drift coefficient of $1$ to the subordinator.
\Cref{fig:convergence_with_drift} shows the numerical convergence rates for $\alpha=0.6$ and $\alpha=0.8$ in this setting, which approach the classical rate of $1/2$.

To illustrate the effectiveness of the truncated EM method for systems with super-linear coefficients, we consider the following numerical example.
	
	\begin{example}\label{ex:2}
		Consider the two-dimensional coupled time-changed SDE system given by
		\begin{equation*}
			\left\{
			\begin{aligned}
				\mathrm{d}X_1(t) &= (X_1^2(t) - 2X_1^5(t))  \mathrm{d}E(t) + X_2^2(t)  \mathrm{d}B_1(E(t)), \\
				\mathrm{d}X_2(t) &= (X_2^2(t) - 2X_2^5(t))  \mathrm{d}E(t) + X_1^2(t)  \mathrm{d}B_2(E(t)), 	
			\end{aligned}
			\right.
		\end{equation*}
		with initial conditions $X_1(0) = 1$ and $X_2(0) = 2$.
	\end{example}
	
	In this example, both the drift and diffusion coefficients exhibit super-linear polynomial growth. 
	Similar to the analysis in \cite{li2025truncated}, it is straightforward to verify that this example satisfies the assumptions of \Cref{thm:TEM_convergence}.

	For the implementation of the truncated EM scheme, we define $\mu(u) = 2u^5$ and $\kappa(\Delta t) = \Delta t^{-\varepsilon}$. 
	We set the truncation parameter to $\varepsilon = 0.01$. 

\begin{figure}[htbp!]
		\centering
		\includegraphics[width=0.48\linewidth]{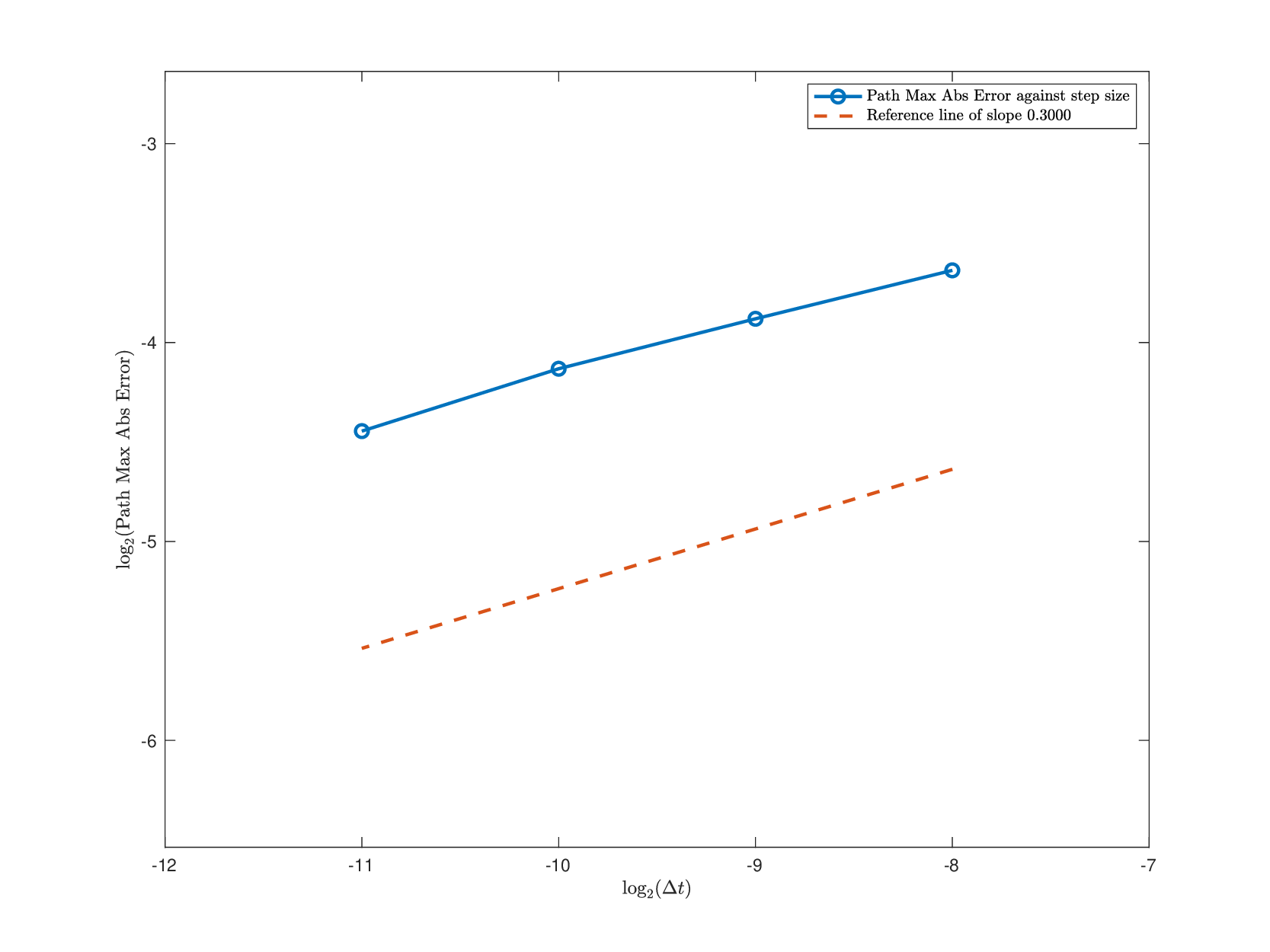}
		\hfill
		\includegraphics[width=0.48\linewidth]{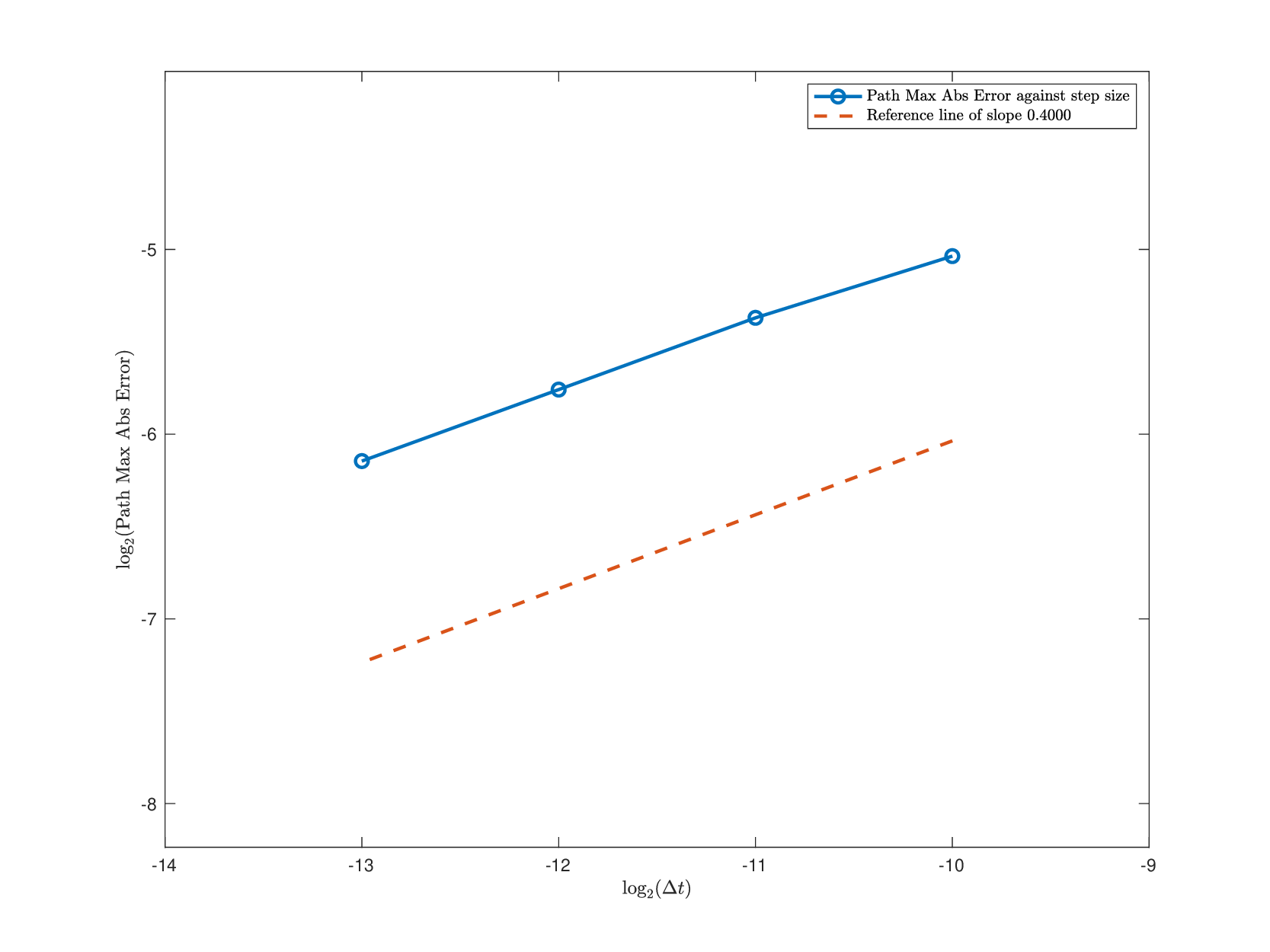}
		\caption{Convergence rates of the truncated EM method for $\alpha=0.6$ (left) and $\alpha=0.8$ (right) in \Cref{ex:2}.}
		\label{fig:TEM_convergence}
\end{figure}

	\Cref{fig:TEM_convergence} displays the $L_1$ errors of the truncated EM method for $\alpha=0.6$ and $\alpha=0.8$. 
	The numerical results show that the empirical convergence rates match the theoretical order of $\alpha/2$, thereby confirming the strong convergence of the proposed scheme for time-changed SDEs with super-linear coefficients.
	
\section*{Declarations}

\bmhead{Funding}
The authors declare that no funding was received from any organization or agency in support of this research.

\bmhead{Data availability}
The data that support the findings of this study are available from the author, upon reasonable request.

\bibliography{reference}

\end{document}